%
%
%
%
%
%
%
%
%
\scrollmode
\magnification=\magstep1
\hoffset=1cm \hsize=12cm

\def\demo#1:{\medskip \it{#1}. \rm}
\def\ni{\noindent}               
\def\ll{\leftline}
\def\cl{\centerline}

%
%
\outer\def\beginsection#1\par{\bigskip
  \message{#1}\centerline{\bf #1}
  \nobreak\smallskip\vskip-\parskip\noindent}

%
%
\outer\def\proclaim#1:#2\par{\medbreak\vskip-\parskip
	{\rm#1.\enspace}{\sl#2}
  \ifdim\lastskip<\medskipamount \removelastskip\penalty55\smallskip\fi}

%
%

\def\R{{\rm I\kern-0.2em R\kern0.2em \kern-0.2em}}
\def\N{{\rm I\kern-0.2em N\kern0.2em \kern-0.2em}}
\def\P{{\rm I\kern-0.2em P\kern0.2em \kern-0.2em}}
\def\B{{\rm I\kern-0.2em B\kern0.2em \kern-0.2em}}
\def\C{{\rm C\kern-.4em {\vrule height1.4ex width.08em depth-.04ex}\;}}
\def\CP{\C\P}
\def\RP{\R\P}

%
%
%
%
\def\cA{{\cal A}}
\def\cB{{\cal B}}
\def\cC{{\cal C}}

\def\cF{{\cal F}}

\def\cH{{\cal H}}

\def\cJ{{\cal J}}

\def\cO{{\cal O}}

\def\cU{{\cal U}}

%
%
%
\def\a{\alpha}
\def\b{\beta}
\def\g{\gamma}
\def\d{\delta}
\def\e{\epsilon}
\def\z{\zeta}

\def\l{\lambda}
\def\r{\rho}

\def\S{{\Sigma}}

%
%
%
%
\def\bar{\overline}              
\def\bs{\backslash}              

\def\di{\partial}                
\def\dibar{\bar\partial}         

%
%
\def\dim{{\rm dim}\,}                    
\def\holo{holomorphic}                   
\def\nbd{neighborhood}                   
\def\psc{pseudoconvex}                   
\def\spsc{strongly\ pseudoconvex}        
\def\ra{real-analytic}                   
\def\spsh{strongly\ plurisubharmonic}
\def\tr{totally real}                    
\def\pc{polynomially convex}             
\def\ss{\subset\!\subset}                
\def\supp{{\rm supp}\,}                  
\def\iff{if and only if}
\def\hvf{holomorphic vector field}

\def\hra{\hookrightarrow}
\def\dist{{\rm dist}}
\def\rank{{\rm rank}}
\def\disc{\triangle}

\def\wt{\widetilde}

\def\ddot#1{{\mathaccent "7F #1}}

\def\begin{\ll{}\vskip 10mm \nopagenumbers}  
\def\pn{\footline={\hss\tenrm\folio\hss}}   

%
%
%
\begin	
\cl{\bf NONCRITICAL HOLOMORPHIC FUNCTIONS}
\cl{\bf ON STEIN MANIFOLDS}
\vskip 5mm
\centerline{by}
\vskip 5mm
\cl{FRANC FORSTNERI\v C}
\bigskip
\it
\cl{IMFM, University of Ljubljana}
\cl{Ljubljana, Slovenia} 
\rm

\vskip 15mm

\cl{\bf 1. Introduction} 
\medskip
\rm

In 1967 Gunning and Narasimhan proved that every open 
Riemann surface admits a \holo\ function without critical points [GN],
thus giving an affirmative answer to a long standing question. Their
proof was an ingenious application of the approximation methods 
of Behnke and Stein.

A complex manifold is called {\it Stein} (after Karl Stein [Ste], 1951) if it is 
biholomorphic to a closed complex submanifold of a complex Euclidean space $\C^N$. 
Open Riemann surfaces are precisely Stein manifolds of complex dimension one.
In this paper we prove the following result.

\proclaim THEOREM I: Every Stein manifold admits a \holo\ function without 
critical points. More precisely, an $n$-dimensional Stein manifold admits 
$[{n+1\over 2}]$ holomorphic functions with pointwise independent 
differentials and this number is maximal for every $n$.

For a more precise statement see Theorems 2.1 and 2.6. An example of  
Forster [Fo1] provides for each $n\in\N$ an  $n$-dimensional Stein manifold 
which does not admit more than $[{n+1 \over 2}]$ holomorphic functions
with independent differentials (Proposition 2.12 below).

The question on the existence of noncritical holomorphic functions on 
a Stein manifold has been open since the 1967 work of 
Gunning and Narasimhan [GN]; it was mentioned in 
Gromov's monograph [Gro3, p.\ 70]. Our proof, which also applies to 
Riemann surfaces, is conceptually different from the one in [GN]. 
It is much easier to construct noncritical smooth real functions 
on smooth open manifolds; see e.g.\ Lemma 1.15 in [God, p.\ 9].

The critical locus of a generically chosen holomorphic function 
on a Stein manifold is discrete. Conversely, we prove that 
for any discrete subset $P$ in a Stein manifold $X$ there exists 
a holomorphic function $f\in \cO(X)$ whose critical locus
equals $P$ (Corollary 2.2).

Recall that a holomorphic map
$f=(f_1,\ldots,f_q)\colon X\to\C^q$ is a {\it submersion} 
if its differential $df_x\colon T_x X\to T_{f(x)}\C^q \simeq \C^q$ 
is surjective for every $x\in X$. Equivalently, the differentials
of its component functions must be linearly independent, i.e., 
$df_1\wedge df_2 \wedge\cdots \wedge df_q\ne 0$. 
Thus the differential of a holomorphic 
submersion $X\to\C^q$ induces a surjective 
complex vector bundle map $TX \to X\times \C^q$ of the tangent
bundle of $X$ onto the trivial bundle of rank $q$ over $X$. 
Our main result is that, for $q<\dim X$, 
this necessary condition for the existence 
of a submersion $X\to\C^q$ is also sufficient.

\proclaim THEOREM II: {\rm (The homotopy principle for holomorphic submersions.)}
If $X$ is a Stein manifold and $1\le q < \dim X$ then every surjective complex 
vector bundle map $TX\to X\times\C^q$ is homotopic 
to the differential of a \holo\ submersion $X\to\C^q$.
\medskip

\pn

The homotopy referred to above belongs to the space 
of surjective complex vector bundle maps  $TX\to X\times\C^q$. 
Theorem II is a holomorphic analogue of the basic 
homotopy principle for submersions of smooth open manifolds 
to real Euclidean spaces, due to A.\ Phillips [Ph1]
and M.\ Gromov [Gro1]. A more precise statement is given by Theorems 2.5 
and 2.6 in Sect.\ 2. We don't know whether the same conclusion 
holds for $q=\dim X>1$ (for open Riemann surfaces see [GN]). 

By using the tools developed in this paper one can also prove 
the following. If $f_0,f_1\colon X\to\C^q$ are \holo\ submersions 
($q<\dim X$) whose differentials $df_0$, $df_1$ are homotopic 
through a family of surjective complex vector bundle maps 
of $TX$ onto the trivial bundle $X\times\C^q$ then there 
exists a homotopy of \holo\ submersions $f_\tau \colon X\to\C^q$ 
$(\tau \in [0,1])$ connecting $f_0$ to $f_1$. 
We shall include the details in a forthcoming 
publication in which we plan to investigate the same problem
for more general target manifolds. 

Theorem I is a corollary of Theorem II and a result of Ramspott [Ra]
to the effect that the cotangent bundle of an $n$-dimensional 
Stein manifold admits $[{n+1\over 2}]$ independent sections, and
these define a surjective complex vector bundle map 
$TX\to X\times \C^{[(n+1)/2]}$. Ramspott's theorem combines the 
Lefschetz theorem [AF] with the standard method of constructing sections 
of fiber bundles over CW complexes by stepwise extension over the skeleta. 
Our proof gives both results simultaneously and does not use
Ramspott's theorem.

We give numerous applications to the existence of nonsingular \holo\ foliations
on Stein manifolds. We prove that every complex vector subbundle 
$E \subset TX$ with trivial quotient $TX/E$ is homotopic to the tangent bundle of 
a \holo\ foliation (Corollary 2.9); the same is true if $TX/E$ admits 
locally constant transition functions (Theorem 7.1). Analogous results for 
smooth foliations on open manifolds were proved by Gromov [Gro1] and 
Phillips ([Ph2], [Ph3], [Ph4]), and on closed manifolds by Thurston ([Th1], [Th2]).
Every $n$-dimensional Stein manifold admits nonsingular holomorphic 
submersion foliations of any dimension $\ge [{n\over 2}]$, and if $X$ has geometric 
dimension $k\le n$ then it admits submersion foliations of any dimension 
$\ge [{k\over 2}]$ (Corollary 2.7). We construct submersion foliations
transverse to certain complex submanifolds of $X$ (Corollaries 2.3 and 2.11)
or containing it as a leaf (Corollaries 2.10 and 7.2).

Our construction depends on three main ingredients developed in this paper.
We postpone the general discussion to Sect.\ 2 and mention at this point 
only the following {\it splitting lemma for biholomorphic maps} (Theorem 4.1):
If $A,B\subset X$ is a  {\it Cartan pair} in a complex manifold $X$ 
then every biholomorphic map $\gamma$ sufficiently uniformly close to the 
identity in a \nbd\ of $A\cap B$ admits a decomposition $\g=\b\circ\a^{-1}$, 
where $\a$ (resp.\ $\b$) is a biholomorphic map close to 
the identity in a \nbd\ of $A$ (resp.\ of $B$). 

This lemma is used to patch a pair of \holo\ submersions $f,g$ to $\C^q$, 
defined in a \nbd\ of $A$ resp.\ $B$, which are sufficiently uniformly close 
in a \nbd\ of $A\cap B$, into a submersion $\wt f$ in a \nbd\ of $A\cup B$. 
The map $\g$ arises as a transition map satisfying $f=g\circ \g$ near $A\cap B$. 
From $\g=\b\circ\a^{-1}$ we obtain $f\circ \a=g\circ \b$ which gives $\wt f$.

Our splitting lemma plays the analogous role in our construction 
of submersions as the {\it Cartan's lemma} (on product splitting 
of holomorphic maps with values in a complex Lie group) does in 
Cartan's theory or in the Oka-Grauert theory. A key difference is that 
our lemma gives a compositional splitting of biholomorphic maps 
and is closer in spirit to Kolmogorov's work  on compositions of 
functions [Ko]. We prove it by a rapidly convergent 
Kolmogorov-Nash-Moser type iteration (Sect.\ 4). 

Our proof of Theorem II breaks down for $q=\dim X>1$ due to a possible 
Picard type obstruction in the approximation problem (Lemma 3.4). 
Hence the following problem remains open.

\medskip
{\it Problem 1.}
Does a parallelizable Stein manifold of dimension $n>1$ holomorphically
immerse in $\C^n$ (i.e., is it a Riemann domain over $\C^n$)~?
\medskip

This well known problem (see [BN, p.\ 18] or [Gro3, p.\ 70]) 
was our main motivation for the present work. To find such an immersion 
it would suffice to obtain an affirmative answer to any of 
the following two problems. 

\medskip 
{\it Problem 2.} 
Let $B$ be an open convex set in $\C^n$ for $n>1$. Is 
every \holo\ immersion ($=$submersion) $B\to \C^n$ 
a uniform limit on compacts of entire immersions $\C^n\to \C^n$~? 
\medskip

The analogous problem for mappings with constant Jacobian may be related to the 
{\it Jacobian problem} for  holomorphic polynomial maps [BN, p.\ 21]. 
The situation is much better understood for biholomorphic maps:
{\it If $f$ is an injective holomorphic map from a convex open set $B\subset \C^n$ 
onto a Runge set $f(B)\subset \C^n$ then $f$ can be approximated uniformly on 
compacts in $B$ by holomorphic automorphisms of $\C^n$} [AL]. No comparable 
result seems to be known for non-injective immersions.

\medskip
{\it Problem 3.} 
Let $f=(f_1,\ldots, f_q) \colon X\to\C^q$ be a \holo\ submersion
for some $q<\dim X$. Given a $(1,0)$-form $\theta_0$ such that 
$df_1\wedge\ldots\wedge df_q\wedge \theta_0\ne 0$ on $X$,
find a homotopy of $(1,0)$-form $\theta_t$ ($t\in [0,1]$) such that
$df_1\wedge\ldots\wedge df_q\wedge \theta_t  \ne 0$ for all $t\in [0,1]$
and $\theta_1=dg$ for some $g\in\cO(X)$. 
(The map $(f,g)\colon X\to\C^{q+1}$ is then a submersion.)
\medskip

Problem 1 has an affirmative answer if one can solve Problem 3 with $q=n-1$.
Explicitly, given a holomorphic submersion $f\colon X^n\to\C^{n-1}$
such that $\ker df$ is a trivial line subbundle of $TX$, find a $g\in \cO(X)$ 
whose restriction to every level set $\{f=c\}$ is noncritical.

%
%
%
%
\beginsection 2. The main results

Let $X$ be a Stein manifold (for their general theory see [GR] and [H\"o2]). 
Denote by $\cO(X)$ the algebra of all \holo\ functions on $X$. 
A compact set $K\subset X$ is said to be 
{\it $\cO(X)$-convex} if for any point $x\in X\bs K$ there 
exists $f\in \cO(X)$ satisfying $|f(x)|> \max_K |f|$. An $\cO(\C^n)$-convex
set is called {\it polynomially convex}. A function 
is holomorphic on a closed subset $K \subset X$ if it is holomorphic 
in some unspecified open \nbd\ of $K$; the set of all such functions
(with the usual identification of functions which agree near $K$)
is denoted $\cO(K)$. We denote by $j^r_x(f)$ the $r$-jet of a function 
$f$ at $x\in X$.  The {\it critical set} of $f\in \cO(X)$ is 
${\rm Crit}(f;X)=\{x\in X\colon df_x=0\}$; a function 
without critical points will be called {\it noncritical}. 
We denote by $|z|$ the Euclidean norm of $z\in \C^n$. 

%
%
\medskip
{\it 1. Functions with prescribed critical locus.} 
Our first main result is

\proclaim THEOREM 2.1: Let $X$ be a Stein manifold, $X_0\subset X$ a closed
complex subvariety of $X$ and $K\subset X$ a compact $\cO(X)$-convex subset.
Let $U\subset X$ be an open set containing $X_0\cup K$ and $f\in \cO(U)$
a holomorphic function with discrete critical set 
$P=\{p_1,p_2,\ldots\}  \subset K\cup X_0$.
For any $\e>0$ and $r, n_1,n_2,\ldots \in \N$ there exists 
an $\wt f\in \cO(X)$ satisfying ${\rm Crit}(\wt f;X)=P$, 
$|f(x)- \wt f(x)|<\e$ for all $x\in K$, $j^r_x(f-\wt f)=0$ 
for all $x\in X_0$, and $j^{n_k}_{p_k}(f-\wt f)=0$
$(k=1,2,\ldots)$. In particular, if $f$ is noncritical on 
$U$ then $\wt f$ is noncritical on $X$.

Theorem 2.1 implies that any noncritical \holo\ function 
on a closed complex submanifold $X_0$ of a Stein manifold $X$ 
extends to a noncritical \holo\ function on $X$. 
Furthermore, there exist noncritical functions satisfying 
the axioms of a Stein manifold ([H\"o2], p.\ 116, Definition 5.1.3.). 
Theorem 2.1 is proved in Sect.\ 5.

The critical locus of a generically chosen \holo\ function on a Stein manifold 
is discrete. Theorem 2.1 implies the following converse.

\proclaim COROLLARY 2.2: Let $P=\{p_1,p_2,p_3,\ldots\}$ be a discrete
set in a Stein manifold $X$ and let $f_k$ be a \holo\ function in a \nbd\ 
of $p_k$ with an isolated critical point at $p_k$ for $k=1,2,\ldots$. 
For any choice of integers $n_k\in \N$ there exists an $f \in \cO(X)$ 
with ${\rm Crit}(f)=P$ such that $f-f_k$ vanishes at least to 
order $n_k$ at $p_k$ for every $k=1,2,\ldots$.

%
%
\medskip
{\it 2. Foliations by complex hypersurfaces.} 
We denote by $TX$ the holomorphic tangent bundle of $X$ and by $T^*X$
its holomorphic cotangent bundle. For the general theory of foliations
we refer to [God].

\proclaim COROLLARY 2.3: Every Stein manifold admits a nonsingular \holo\ 
foliation by closed complex hypersurfaces; in addition such a foliation 
may be chosen to be transverse to a given closed complex submanifold.

\demo Proof:
A closed complex submanifold $V$ of a Stein manifold $X$ is itself 
Stein and hence admits a noncritical function $f\in \cO(V)$ by Theorem 2.1. 
By Cartan's theorem $f$ extends 
to a holomorphic function on $X$. Since the extension remains
noncritical on $X_0$, Theorem 2.1 gives a noncritical function 
$\wt f\in\cO(X)$ such that $\wt f|_V=f$. The family of connected 
components of the levels sets $\{\wt f=c\}$ ($c\in \C$) is a 
foliation of $X$ by closed complex hypersurfaces transverse to $V$. 

\proclaim COROLLARY 2.4: 
If $V$ is a smooth closed complex hypersurface with trivial
normal bundle in a Stein manifold $X$ then $V$ is a union of 
leaves in a nonsingular holomorphic foliation of $X$ by closed complex 
hypersurfaces. This holds in particular if $H^2(V;Z)=0$, or if $X=\C^n$.
Any smooth connected complex curve in a Stein surface 
is a leaf in a nonsingular \holo\ foliation.

{\it Proof.} 
Choose a holomorphic trivialization of the normal 
bundle $N=TX|_V/TV \simeq V\times\C$. The projection 
$h\colon N\to\C$ on the second factor is a noncritical 
\holo\ function on $N$ and $N_0=\{h=0\}$ is  
the zero section of $N$. The Docquier-Grauert theorem
[DG] (see also Theorem 8 in [GR, p.\ 257]) 
gives an open \nbd\ $\Omega\subset X$ of $V$ and an injective 
holomorphic map $\phi \colon \Omega\to N$ with $\phi(V)=N_0$. 
Then $f=h\circ \phi$ is a noncritical function
on $\Omega$ with $\{f=0\}=V$. Applying Theorem 2.1 (with $X_0=V$) we 
obtain a noncritical function $\wt f\in\cO(X)$ which vanishes on $V$,
and the foliation $\{\wt f=c\}$ clearly satisfies Corollary 2.4. 
The second statement follows from the isomorphism 
$Pic(V)= H^1(V;\cO^*) \simeq H^2(V;Z)$; the latter group 
vanishes if $V$ is an open Riemann surface. Since
every divisor on $\C^n$ is a principal divisor, the normal
bundle of any complex hypersurface $V\subset \C^n$ is trivial.
\medskip

%
%
%
%
{\it 3. Holomorphic submersions and foliations.}
We now consider the existence of holomorphic submersions 
$f=(f_1,\ldots,f_q)\colon X\to\C^q$ for $q\le n=\dim X$. 
The components of $f$ are noncritical functions
with pointwise independent differentials, i.e.,
$df_1\wedge \cdots\wedge df_q \ne 0$ on $X$. Hence 
an obvious necessary condition is that there exists a 
$q$-tuple $\theta=(\theta_1,\ldots,\theta_q)$ of continuous differential 
$(1,0)$-forms on $X$  satisfying 
$\theta_1\wedge\cdots\wedge \theta_q|_x \ne 0$ for all $x\in X$. 
Any such ordered $q$-tuple will be called a {\it q-coframe\/} on $X$. 
We may view $\theta$ as a {\it complex vector bundle epimorphism} 
$\theta\colon TX\to X\times\C^q$ of the tangent bundle $TX$ onto the 
trivial bundle of rank $q$ over $X$. Clearly we may speak of 
{\it holomorphic $q$-coframes, homotopies of $q$-coframes}, etc.
If $\theta_j=df_j$ for some $f_j\in \cO(X)$ $(j=1,\ldots,q)$
we shall write $\theta=df$ and call $\theta$ {\it exact holomorphic}. 
The following two theorems are our main results;
they are proved in Sect.\ 6.

\proclaim THEOREM 2.5: 
Let $X$ be a Stein manifold and $1\le q<\dim X$. For every $q$-coframe 
$\theta^0$ on $X$ there exists a homotopy of $q$-coframes 
$\theta^t$ $(t\in [0,1])$ such that $\theta^1=df$ where 
$f\colon X\to \C^q$ is a holomorphic submersion. Furthermore, 
if $X_0, K\subset X$ are as in Theorem 2.1, $r\in \N$, $\e>0$, 
and if we assume that $\theta^0=df^0$ is exact holomorphic in an open set 
$U\supset X_0\cup K$, the homotopy may be chosen such that 
$\theta^t=df^t$ is exact holomorphic in a \nbd\ of $X_0\cup K$ for every 
$t\in[0,1]$, $f^t-f^0$ vanishes to order $r$ on $X_0$, and $|f^t-f^0|<\e$ on $K$.

Theorem 2.5 also holds for $q=\dim X=1$ and is due in this case 
to Gunning and Nara\-sim\-han who proved that for every nonvanishing holomorphic
one-form $\theta$ on an open Riemann surface there exists 
a holomorphic function $w$ such that $e^w\theta=df$ is exact
holomorphic [GN, p.\ 107]. The homotopy $\theta^t=e^{tw}\theta$ 
consisting of nonvanishing one-forms connects $\theta^0=\theta$ 
to $\theta^1=df$. We do not know whether Theorem 2.5 holds for $q=\dim X>1$.

We state separately the case when the necessary condition 
on the existence of a $q$-coframe is automatically fulfilled
due to topological reasons.
Recall that any Morse critical point of a \spsh\ function on an 
$n$-dimensional complex manifold has Morse index at most $n$ [AF].
If $X$ admits a \spsh\ Morse exhaustion function $\rho\colon X\to \R$ 
all of whose critical points have index $\le k$ (and $k$ is 
minimal such), we say that $X$ has {\it geometric dimension $k$};
such $X$ is homotopically equivalent to a $k$-dimensional 
CW-complex [AF].

\proclaim THEOREM 2.6: Let $\rho\colon X\to \R$ be a \spsh\ Morse exhaustion 
function on an $n$-dimensional Stein manifold $X$.  
Assume that $c$ is a regular value of $\rho$ and every critical point
of $\rho$ in $\{x\in X\colon \rho(x) >c\}$ has Morse index $\le k$. 
If $q\le q(k,n):= \min\{n-[{k\over 2}],n-1\}$ then every \holo\ submersion 
$f\colon \{x\in X\colon \rho(x) <c\} \to\C^q$ can be approximated 
uniformly on compacts by \holo\ submersions $\wt f\colon X\to \C^q$. 
Every $n$-dimensional Stein manifold $X$ admits a holomorphic 
submersion to $\C^{[(n+1)/2]}$; if $X$ has geometric dimension 
$k$ then it admits a holomorphic submersion to $\C^{q(k,n)}$.

\smallskip
Proposition 2.12 below shows that the submersion dimension in 
Theorem 2.6 is optimal for every $n$. Theorem 2.6 immediately gives 

\proclaim COROLLARY 2.7: Every Stein manifold $X$ of geometric dimension 
$k$ admits nonsingular \holo\ foliations of any dimension $\ge [{k\over 2}]$.
If $X$ is parallelizable, it admits a \holo\ submersion $X\to \C^{n-1}$ 
$(n=\dim X)$ and nonsingular holomorphic foliations of any dimension $\ge 1$.

The foliations in Corollary 2.7 are given by submersions to Euclidean spaces; 
hence all leaves are topologically closed and the normal bundle 
is trivial.

\smallskip
{\it Remark.} 
The Oka-Grauert principle applies to $q$-coframes on a Stein manifold
and shows that any $q$-coframe is homotopic to a \holo\ $q$-coframe,  
and any homotopy between a pair of holomorphic $q$-coframes can be deformed to 
a homotopy consisting of holomorphic $q$-coframes. This is seen 
by viewing $q$-coframes as a sections of the holomorphic fiber 
bundle $V^q(T^*X)\to X$ whose fiber $V^q_x$  is the Stiefel variety 
of all ordered $q$-tuples of $\C$-independent elements in $T^*_x X$. Since the 
Lie group $GL_n(\C)$ $(n=\dim X)$ acts transitively on $V_x^q$, 
the Oka-Grauert principle [Gra] applies to sections $X\to V^q(T^* X)$.

%
%
\medskip
{\it 4. Existence of homotopies to integrable subbundles.} 
The components of a $q$-coframe on $X$ are linearly independent
sections of $T^*X$ which therefore span a  trivial complex subbundle 
of rank $q$ in $T^*X$. Conversely, every trivial rank $q$ 
subbundle $\Theta \subset T^*X$ is spanned by (the components 
of) a  $q$-coframe. Different $q$-coframes $\theta,\theta'$ spanning 
the same subbundle of $T^*X$ are related by 
$\theta' = \theta \cdotp A$ for some  $A\colon X\to GL_q(\C)$. 
A homotopy of $q$-coframes induces a homotopy of the associated subbundles 
of $T^* X$. Hence Theorem 2.5 implies

\proclaim COROLLARY 2.8: 
Let $X$ be a Stein manifold. Every trivial complex subbundle $\Theta\subset T^*X$ 
of rank $q<\dim X$ is homotopic to a subbundle generated 
by independent \holo\ differentials $df_1,\ldots,df_q$. If $\Theta$ 
is \holo\ then the homotopy can be chosen through \holo\ subbundles
of $T^* X$.

The last statement follows from the Oka-Grauert principle [Gra].
Corollary 2.8 admits the following dual formulation in terms of 
subbundles of $TX$ (for a generalization see Theorem 7.1).

\proclaim COROLLARY 2.9: 
Let $X$ be a Stein manifold of dimension $n$. Every complex subbundle 
$E\subset TX$ of rank $k\ge 1$ with trivial quotient bundle $TX/E$ is homotopic 
to an integrable \holo\ subbundle $\ker df \subset TX$, where
$f\colon X\to \C^{n-k}$ is a \holo\ submersion. If $E$ is holomorphic 
then the homotopy may be chosen through \holo\ subbundles.

\demo Proof:
The complex subbundle $\Theta =E^\perp \subset T^*X$ with fibers 
$\Theta_x =\{\l \in T^*_x X\colon \l(v)=0\ {\rm for\ all}\ v\in E_x\}$
(the {\it conormal bundle} of $E$) satisfies 
$\Theta\simeq (TX/E)^*$ and hence is trivial. Corollary 2.8 gives a homotopy 
$\Theta^t\subset T^*X$ $(t\in [0,1])$ from $\Theta^0=\Theta$ to a subbundle 
$\Theta^1 \subset T^*X$ spanned by $n-k$ independent 
holomorphic differentials $df_1,\ldots, df_{n-k}$. The homotopy 
$E^t=(\Theta^t)^\perp \subset TX$ satisfies Corollary 2.9.
The last statement follows from the Oka-Grauert principle.
\medskip

We conclude with a couple of results on the existence of 
submersion foliations which either contain a given submanifold as a leaf,
or else are transverse to it. Both depend on Theorem 2.5 and are proved
in Sect.\ 6.

\proclaim COROLLARY 2.10:
Let $X$ be an $n$-dimensional Stein manifold and $V \subset X$ a closed complex 
submanifold. If $TX$ admits a trivial complex subbundle 
$N$ satisfying $TX|_V= TV\oplus N|_V$ then there is a holomorphic submersion 
$f\colon X\to\C^q$ $(q=n -\dim V)$ such that $V$ is a union of connected components of 
the fiber $f^{-1}(0)$. If $\dim V\ge [{n \over 2}]$ then the above conclusion holds 
provided $V$ has a trivial normal bundle in $X$.

\proclaim COROLLARY 2.11: Let $X$ be a Stein manifold, 
$\iota\colon V\hra  X$ a closed complex submanifold, and 
$f=(f_1,\ldots,f_q) \colon V\to \C^q$ a holomorphic submersion. 
If there is a $q$-coframe $\theta=(\theta_1,\ldots,\theta_q)$ 
on $X$ satisfying $\iota^*\theta_j=df_j$ $(j=1,\ldots,q)$ 
then there exists a \holo\ submersion $F\colon X\to \C^q$ 
with $F|_V=f$. Such $F$ always exists if 
$q\le [{n+1\over 2}]$, where $n=\dim X$.

%
%
\medskip
{\it 5. An example.} The following example shows that the submersion dimension
in Theorem 2.6 is maximal for every $n$.

\proclaim PROPOSITION 2.12: Set
$Y = \{[x\colon y\colon z]\in \CP^2\colon x^2+y^2+z^2\ne 0\}$ and
$$
	X= \cases { Y^m, & if $n=2m$;\cr
                          Y^m\times\C, & if $n=2m+1$. \cr}                        
$$
Then $X$ is an $n$-dimensional Stein manifold which does 
not admit a \holo\ submersion to $\C^{[(n+1)/2]+1}$.

\demo Proof: 
These manifolds were considered by Forster [Fo1, p.\ 714], [Fo2, Proposition 3].
He showed that $Y$ is a Stein surface which admits a strong 
deformation retraction onto the real projective plane 
$M=\{[x\colon y\colon z] \colon x,y,z\in \R\} \simeq \RP^2$ contained in $Y$ as 
a \tr\ submanifold. Thus $Y$ is a complexified $\RP^2$ and $X$ is homotopic to 
$(\RP^2)^m$. Using the fact that $TY|_M \simeq TM\oplus TM$  (as real bundles)
Forster proved that the Stiefel-Whitney class $w_{2m}(TX)$ is the nonzero 
element of the group $H^{2m}(X;Z_2) = H^{2m}((\RP^2)^m; Z_2)=Z_2$
and consequently the Chern class $c_m(TX)$ is the nonzero element 
of $H^{2m}(X;Z)=Z_2$. Hence $c_m(T^*X)=(-1)^m c_m(TX)\ne 0$ [MS, p.\ 168] 
which implies that $T^*X$ does not contain a trivial complex subbundle of 
rank $n-m+1 = [{n+1 \over 2}]+1$. 
(Proof: if $T^*X=E\oplus E'$ where $E'$ is trivial then $0\ne c_m(T^*X)=c_m(E)$ 
[MS, Lemma 14.3] which means that $\rank E \ge m$ and consequently 
$\rank E' \le n-m$.) Hence there exists no submersion $X\to \C^{[(n+1)/2]+1}$.
(The only essential property of $X$ is that the Chern class of $TX$ 
of order $[{n\over 2}]$ does not vanish.)

\medskip
Recall that holomorphic immersions of a Stein manifold $X$ into 
Euclidean spaces of dimension $N>\dim X$ satisfy the 
following homotopy principle (Eliashberg and Gromov [Gro3, pp.\ 65-75]): 
{\it Every injective complex vector bundle map $TX\to X\times \C^N$ 
is homotopic to the differential of a \holo\ immersion $X\to\C^N$.}
In particular, every $n$-dimensional Stein manifold admits a \holo\ 
immersion in $\C^{[3n/2]}$, and the manifold $X$ in Proposition 2.12 does not 
immerse in $\C^{[3n/2]-1}$ [Fo2, p.\ 183]. A comparison with Theorem 2.6 
shows that the {\it submersion dimension} $q(n)$, resp.\ the 
{\it immersion dimension} $N(n)$, are symmetric with respect to $n=\dim X$:
$$ 
	q(n)= \biggl[{n+1\over 2}\biggr]= n- \biggl[{n\over 2}\biggr],
	\quad 	N(n) = n+ \biggl[{n\over 2}\biggr]. 
$$ 
If $X$ has geometric dimension at most $k$ and $k\ge 2$ 
then $X$ admits a submersion to $\C^{n-[k/2]}$ and immersion in $\C^{n+[k/2]}$,
and both bounds are sharp (an example is the manifold
$Y^{[k/2]}\times \C^{n-2[k/2]}$ where $Y$ is as in Proposition 2.12).

%
%
\medskip
{\it 6. Remarks on parallelizable Stein manifolds.}
By Grauert [Gra] the tangent bundle of a Stein manifold $X$ is holomorphically 
trivial if and only if it is topologically trivial (as a complex vector bundle). 
The question whether every such manifold immerses in $\C^n$ with $n=\dim X$ 
remains open for $n>1$. Every closed complex submanifold $X\subset \C^N$ 
with trivial normal bundle is parallelizable [Fo1, p.\ 712]. (Triviality of 
the normal bundle is equivalent to $X$ being a \holo\ complete intersection 
in some open \nbd.) In particular, every closed complex hypersurface in 
$\C^{n+1}$ is parallelizable [Fo1, Corollary 2] but it is unknown whether 
these immerse into $\C^n$. J.\ J.\ Loeb [BN, p.\ 19] found explicit holomorphic 
immersions $X\to \C^n$ of algebraic hypersurfaces $X\subset \C^{n+1}$
of the following type:
$$
	X=\{(z_0,z_1,\ldots,z_k) \colon z_0^d + P_1(z_1)+\ldots + P_k(z_k)=1 \} 
	\subset\C^{n+1},
$$
where $z_0\in\C$, $z_j\in \C^{n_j}$, $P_j$ is a homogeneous polynomial of some 
degree $d_j$ on $\C^{n_j}$ for every $j=1,\ldots,k$, and $n_1+\ldots n_k=n$.
These manifolds are even algebraically parallelizable but do not admit 
algebraic immersions to $\C^n$. An example of this type is the complex $n$-sphere 
$\Sigma^n= \{z\in \C^{n+1}\colon \sum z_j^2 = 1\}$.

In another direction, Y.\ Nishimura [N] found explicit holomorphic immersions
$\CP^2 \bs C\to \C^2$ where $C$ is an irreducible cuspidal cubic in 
$\CP^2$. Further examples and remarks on parallelizable Stein manifolds 
can be found in [Fo1].  

Unlike $\Sigma^n$, the real $n$-sphere $S^n=\Sigma^n\cap \R^{n+1}$ (which is a 
maximal totally real submanifold of $\S^n$) is parallelizable only for $n=1,3,7$. 
By Thurston ([Th1], [Th2]) $S^3$ and $S^7$ admit $\cC^\infty$ foliations of 
all dimensions. However, a simply connected closed \ra\ manifold 
(such as $S^n$ for $n>1$) does not admit  any \ra\ foliations of codimension one 
(Haefliger [Ha1]).

%
%
\medskip
{\it 7. Comparison with smooth immersions and submersions.}
The homotopy classification of smooth immersions $X\to \R^q$ was 
discovered by Smale [Sm] and Hirsch ([Hi1], [Hi2]). 
Subsequently analogous results were proved for submersions 
(Phillips [Ph1] and Gromov [Gro1]), $k$-mersions (Feit [Fe]), 
and maps of constant rank [Ph5]. 
Gromov's monograph [Gro3] offers a comprehensive survey;
see also the more recent monographs [Sp] and [EM]. Our Theorem 2.5 is a 
holomorphic analogue of the basic homotopy principle for smooth submersions 
$X\to\R^q$ which holds for all $1\le q\le \dim_{\R} X$ provided that $X$ 
is a smooth open manifold (see [Hi2], [Ph1], [Gro1], and [Ha2]). 

The differential relation controlling smooth immersions of positive 
codimension is {\it ample in the coordinate directions}, and the corresponding 
homotopy principle follows from the {\it convex integration lemma} 
of M.\ Gromov (see the discussion and references in Subsect.\ 6.2 below). 
The smooth submersion relation is not ample in the coordinate
directions (Example 2 in [EM, p.\ 168]), and the homotopy principle 
for smooth submersions is obtained by exploiting the invariance 
of the submersion condition under local diffeomorphisms and reducing 
the problem to a subpolyhedron in the given manifold.
On the other hand, we shall see that the complex (holomorphic) submersion relation 
is ample in the coordinate directions on any totally real submanifold,
and this is exploited to obtain a maximal rank extension of 
the map across a totally real handle (Lemma 6.5). The invariance under 
local biholomorphisms is also strongly exploited in the approximation 
and patching of submersions.

The homotopy principle type results on Stein manifolds are traditionally 
referred to as (instances of) the {\it Oka principle}; see the 
recent survey [F3].

%
%
\medskip
{\it 8. Outline of proof of the main theorems.}
Our construction relies on three main ingredients developed in this paper. 

The first one is a new technique for approximating a noncritical holomorphic function 
$f$ on a compact polynomially convex subset $K\subset \C^n$  
by entire noncritical functions (Sect.\ 3). We exploit the invariance of the 
maximal rank condition under biholomorphisms. Choose a preliminary 
approximation of $f$ on $K$ by a holomorphic polynomial $h$ with finite critical 
set $\Sigma={\rm Crit}(h)$ disjoint from $K$. When $n>1$, the main step is 
to find an injective holomorphic map $\phi \colon \C^n \to\C^n\bs \Sigma$ 
(a Fatou-Bieberbach map) which is close to the identity map on $K$ and whose 
range avoids $\Sigma$. Such $\phi$ can be obtained as a limit of 
holomorphic automorphisms of $\C^n$ using methods developed by Andersen and Lempert 
([A], [AL]) and  Rosay and the author ([FR], [F1], [F2]). Then $\wt f= h\circ \phi$ is 
noncritical on $\C^n$ and approximates $f$ uniformly on $K$. For $n=1$ we give 
a different proof using Mergelyan's theorem. Similar methods are developed 
for submersions $\C^n\to \C^q$ for $q<n$.

The second ingredient concerns patching of holomorphic submersions.
Let $A,B\subset X$ be compact sets in a complex manifold $X$ 
such that $A\cup B$ has a basis of Stein \nbd s and 
$\bar{A\bs B}\cap \bar{B\bs A}=\emptyset$. For any biholomorphic 
($=$injective holomorphic) map $\gamma \colon V\to X$, sufficiently close to the 
identity map in a \nbd\ $V\subset X$ of $C=A\cap B$, we obtain a compositional splitting 
$\gamma=\beta\circ\alpha^{-1}$, where $\a$ (resp.\ $\b$) is a biholomorphic 
map close to the identity in a \nbd\ of $A$ resp.\ of $B$ (Theorem 4.1). 
If $f$ (resp.\ $g$) is a submersion to $\C^q$ in a \nbd\ of $A$ (resp.\ $B$) 
and $g$ is sufficiently uniformly close to $f$ in a \nbd\ of $C$ then 
$f=g\circ \gamma$ for a biholomorphic map $\gamma$ close to the identity;
splitting $\gamma=\beta\circ\alpha^{-1}$ as above we obtain $f\circ \alpha=g\circ \beta$ 
near $C$; this gives a submersion $\wt f$ in a \nbd\ of $A\cup B$ which 
approximates $f$ on $A$. 

\smallskip
{\it Remark.} 
The standard $\dibar$-method for patching $f$ and $g$ would be to take 
$h=f + \chi(g-f)$ and $\wt f= h - T(\dibar h)$, where $T$ is a bounded solution operator 
to the $\dibar$-equation in a \nbd\ of $A\cup B$ and $\chi$ 
is a smooth function which equals zero in a \nbd\ of $\bar {A\bs B}$ 
and one in a \nbd\ of $\bar{B\bs A}$. Since $\dibar h= (g-f)\dibar \chi$, 
the correction term $T(\dibar h)$ is controlled by $|f-g|$ and hence $\wt f$ 
is noncritical in a \nbd\ of $A$ provided that $|g-f|$ is sufficiently 
small in a \nbd\ of $A\cap B$. However, to insure that $\wt f$ is also noncritical 
in a \nbd\ of $B$ we would need the pointwise estimate $|d(T(\dibar h))| < |dg|$. 
Since we obtain $g$ by a Runge approximation of $f$ on $A\cap B$ (and we have 
no control on its differential on $B\bs A$), such an estimate is impossible.
\smallskip

In the construction of submersions $X\to\C^q$ for $q>1$ another nontrivial problem
is the crossing of the critical levels of a \spsh\ Morse exhaustion function 
$\rho$ on $X$. We combine three ingredients (Sect.\ 6):
\smallskip
\item{--} a {\it convex integration lemma} of Gromov, or Thom's  
{\it jet transversality theorem} when $q\le [{n+1\over 2}]$,
to obtain a smooth extension across a handle,
\item{--} \holo\ approximation on certain handlebodies, and 
\item{--} the construction of an increasing family of smooth  
\spsc\ \nbd s of a handlebody, passing over the critical level
of $\rho$.
\smallskip

We globalize the construction using the `bumping method' similar to the 
one in [Gro4], [HL3], [FP1], [FP2], [FP3]. We exhaust $X$ by an increasing sequence 
$A_0\subset A_1\subset A_2\subset \cdots\subset \cup_{k=1}^\infty A_k = X$
of compact $\cO(X)$-convex sets such that the initial function 
(or submersion) $f=f_0$ is defined on $A_0$, and for each $k\ge 0$ 
we have $A_{k+1}=A_k\cup B_k$ where $(A_k,B_k)$ is a 
{\it special Cartan pair}. This enables us to approximate a noncriticial 
function $f_k$ on $A_k$ by a noncritical function $f_{k+1}$ on 
$A_{k+1}$. The limit $\wt f=\lim_{k\to\infty} f_k$ 
is a noncritical function on $X$. The details are given
in Sect.\ 5 for functions and in Sect.\ 6 for submersions.

In Section 7 we construct holomorphic sections transverse to certain 
\holo\ foliations, thus generalizing Corollaries 2.9 and 2.10.

%
%
%
%
\beginsection 3. Approximation of noncritical functions and submersions

This section uses  the Anders\'en-Lempert 
theory of holomorphic automorphisms of $\C^n$ ([A], [AL]) as developed further 
in ([FR], [F1], [F2]). The following is one of the main steps in our construction 
of noncritical holomorphic functions.

\proclaim THEOREM 3.1: Let $K$ be a compact \pc\ subset of $\C^n$.
Let $f$ be a \holo\ function in an open set $U\supset K$ satisfying $df\ne 0$. 
Given $\e>0$ there exists a $g\in\cO(\C^n)$ 
satisfying $dg\ne 0$ on $\C^n$ and $\sup_K |f-g| <\e$.

\demo Proof: 
Choose a compact \pc\ set $L\subset U$ with smooth boundary and containing 
$K$ in the interior. Such $L$ may be obtained as a regular sublevel set of a 
\spsh\ exhaustion function on $\C^n$ which is negative on $K$ and positive on 
$\C^n\bs U$ ([H\"o2], Theorem 2.6.11.).

Consider first the case $n=1$. Since $L \subset \C$ is smoothly bounded and 
\pc, it is a union $L=\cup_{j=1}^m L_j$ of finitely many compact, connected 
and simply connected sets $L_j$. Since $f'(z)\ne 0$ for $z\in U$, 
there is a \holo\ function $h$ in a \nbd\ of $L$ such that 
$f'(z)=e^{h(z)}$ for each $z$. 

For every $j=2,\ldots,m$ we connect $L_1$ to $L_j$ by a simple smooth arc 
$C_j$ contained in $\C\bs L$ except for its endpoints
$a_j\in L_1$, $b_j\in L_j$. Furthermore we choose the arcs $C_j$ to be 
pairwise disjoint. The sets $S:=L\cup C_2\cup\cdots\cup C_m$ and 
$\C\bs S$ are connected, and $h$ can be extended to a smooth function 
on $C_j$ satisfying $\int_{C_j}  e^{h(\z)}\, d\z = f(b_j)-f(a_j)$ for 
$j=2,\ldots, m$ (where $C_j$ is oriented from $a_j$ to $b_j$).
By Mergelyan's theorem we can approximate $h$ uniformly on $S$ as close as 
desired by a \holo\ polynomial $\wt h$. Choose a point $a\in L_1$ and define 
$g(z)=f(a)+\int_{a}^z e^{\wt h(\zeta)}\, d\zeta$. The integral 
does not depend on the choice of the path and hence $g$ is an entire 
function on $\C$ satisfying $g'(z)=e^{\wt h(z)} \ne 0$ for each $z\in \C$. 
If $z\in L$, we can choose the path of integration from $a$ to $z$  
entirely contained in $S$ and with length bounded from above 
independently of $z$. (If $z\in L_j$ for $j>1$, we include the arc $C_j$ 
in the path of integration.) It follows that $g$ approximates $f$ uniformly 
on $L$. This completes the proof for $n=1$.

Assume now $n\ge 2$. Since $L$ is \pc, there exists for any $\e>0$ 
a \holo\ polynomial $h$ on $\C^n$ satisfying $\sup_L |f-h|<\e/2$. 
If $\e$ is chosen sufficiently small then $dh\ne 0$ on $K$. For a
generic choice of $h$ its critical set 
$\Sigma=\{z\in\C^n \colon dh_z=0\} \subset \C^n\bs K$ is finite
(since it is given by $n$ polynomial equations 
$\di h/\di z_j=0$, $j=1,\ldots, n$). 
To complete the proof we need the following.

\proclaim PROPOSITION 3.2: Let $K$ be a compact \pc\ subset of $\C^n$
$(n\ge 2)$. Given a finite set $\Sigma \subset \C^n\bs K$ and $\d>0$ 
there is a biholomorphic map $\phi$ of $\C^n$ onto a subset 
$\Omega\subset \C^n\bs \Sigma$ such that $|\phi(z)-z|<\d$ for all $z\in K$.

Recall that a biholomorphic of $\C^n$ onto a proper subset of $\C^n$ 
is called a {\it Fatou-Bieberbach map}. Thus $\phi$ is a Fatou-Bieberbach 
map whose restriction to $K$ is close to the identity map and whose 
range avoids $\Sigma$.

Assume for a moment that Proposition 3.2 holds. Let $c=\sup_{z\in L}|dh_z|$. 
Choose $\d <\min\{\dist(K,\C^n\bs L), \e/2c\}$.
Let $g=h\circ\phi \in\cO(\C^n)$ where $\phi$ is furnished by
Proposition 3.2. Then $dg_z=dh_{\phi(z)}\cdotp d\phi_z \ne 0$ for 
every $z\in\C^n$ (since $\phi(z)\in\Omega \subset\C^n\bs \Sigma$ and 
$dh\ne 0$ on $\C^n\bs \Sigma$). For every $z\in K$ we have
$$
	|g(z)-h(z)| = |h(\phi(z))-h(z)| \le c|\phi(z)-z| < c\d < \e/2 
$$
and hence $|g(z)-f(z)|<\e$. This proves Theorem 3.1.

\demo Proof of Proposition 3.2:
Choose $\e\in (0,1)$. Let $B$ denote the closed unit ball centered at the origin 
in $\C^n$ and $rB$ its dilation by $r>0$. 
Choose a compact set $L\subset \C^n\bs \Sigma$ containing 
$K$ in its interior. Let $r_1>1$ be chosen such that $L\subset (r_1-1)B$.
Set $r_k=r_1+k-1$ and $\e_k=2^{-k-1}\e$ for $k=1,2,3,\ldots$.   

Consider the \holo\ flow on a \nbd\ of $L\cup \Sigma$ in $\C^n$ which
rests near $L$ and moves the finite set $\Sigma$ out of the ball $r_1 B$.
Since the trace of this flow is \pc, the time-one map can be 
approximated uniformly on $L$ by \holo\ automorphisms of $\C^n$
according to Theorem 1.1 in [FR]. This gives a holomorphic automorphism 
$\psi_1$ of $\C^n$ satisfying $|\psi_1(z)-z| <\e_1$ for $z\in L$ and 
$\psi_1(\Sigma) \cap r_1B =\emptyset$. (That is, we pushed $\Sigma$ 
out of the ball $r_1 B$ by a holomorphic automorphism of $\C^n$ which is 
$\e_1$-close to the identity map on $L$.) 

Set $\Sigma_1=\psi_1(\Sigma)$. By the same argument there is 
an automorphism $\psi_2$ of $\C^n$ satisfying $|\psi_2(z)-z|<\e_2$ for 
$z\in r_1 B$ and  $\psi_2(\Sigma_1) \cap r_2B =\emptyset$. 

Continuing inductively we obtain a sequence of automorphisms $\psi_k$ of 
$\C^n$ such that $|\psi_k(z)-z|<\e_k$ on $r_{k-1} B$ and 
$\psi_k(\Sigma_{k-1}) \cap r_k B=\emptyset$ for each $k=1,2,\ldots$.
By Proposition 5.1 in [F2] (which is entirely elementary)
the sequence of compositions   
$\psi_k\circ\psi_{k-1}\circ\cdots\circ\psi_1$ converges as 
$k\to\infty$ to a biholomorphic map $\psi\colon \Omega\to\C^n$ 
from an open set $\Omega\subset \C^n$ onto $\C^n$.
By construction we have $L\subset \Omega \subset \C^n\bs \Sigma$
and $|\psi(z)-z|<\e$ for $z\in L$. The inverse map
$\phi=\psi^{-1}\colon \C^n\to\Omega$ is biholomorphic onto 
$\Omega \subset \C^n\bs \Sigma$ and is uniformly close to the identity on $K$. 
Choosing $\e$ sufficiently small we can insure that $|\phi(z)-z|<\d$ for $z\in K$.
This proves Proposition 3.2.
\medskip

To construct holomorphic submersions $X\to\C^q$ for $1<q<n=\dim X$ we need 
a suitable analogue of Theorem 3.1 for submersions $f\colon U\to \C^q$,
where $U$ is an open set in $\C^n$ containing a given compact 
\pc\ set $K\subset \C^n$. The initial approximation of $f$ gives a polynomial
map $h\colon\C^n\to\C^q$ for which $\Sigma:=\{z\in\C^n\colon \rank\, dh_z <q\}$
is an algebraic subvariety of $\C^n$ of complex dimension $q-1$ (which 
is at most $n-2$). To conclude the proof as above we would need a
Fatou-Bieberbach map whose range  contains $K$ but omits $\Sigma$. 
Unfortunately we have been unable to construct such a map, 
and we even have some doubts about its existence due to the possible 
linking of $K$ and $\Sigma$. Instead we prove a result of this kind only 
for very special pairs $(K,\Sigma)$ which suffices for the application at hand.

\proclaim PROPOSITION 3.3: 
Let $x=(z,w)$ be complex coordinates on $\C^n=\C^r\times \C^s$. 
Let $D\subset \C^r$ and $K\subset \C^n$ be compact \pc\ sets
such that $D\times\{0\}^s \subset K\subset D\times \C^s$ and 
each fiber $K_z=\{w\in \C^s\colon (z,w)\in K\}$ $(z\in D)$ is convex. 
Assume that $s\ge 2$ and $q \le r+1$. Let $f\colon U\to\C^q$ 
be a \holo\ submersion in an open set $U\subset\C^n$ containing $K$.
Given $\e>0$ and a compact set $L\subset D\times \C^s$ containing $K$, there 
exists a \holo\ submersion $g \colon V \to \C^q$ in an open 
set $V\supset L$ satisfying $\sup_K |f-g| <\e$.

\medskip {\it  Remark.}
Proposition 3.3 is only used in the proof of Proposition 6.1 (Sect.\ 6)
with $r=n-2$, $s=2$; hence the only condition on $q$ is $q\le n-1$.

\demo Proof: 
Denote by $\pi\colon\C^n\to\C^r$ the projection $\pi(z,w)=z$. 
We can approximate $f$ uniformly on a \nbd\ of $K$ by a 
polynomial map $h=(h_1,\ldots,h_q)\colon \C^n\to\C^q$. A generic choice of $h$ 
insures that the set $\Sigma := \{x\in\C^n\colon \rank\, dh_x <q\}$ is an algebraic 
subvariety of dimension $q-1 \le r$ which does not intersect $K$ and the projection 
$\pi|_\Sigma \colon \Sigma\to\C^r$ is proper. (This follows from the jet 
transversality theorem: $\Sigma$ is the common zero set of all maximal minors 
of the complex $q\times n$ matrix $\pmatrix{\di h_j/\di x_l}$; at each point 
at least $n-(q-1)$ of these equations are independent. Hence for a generic choice 
of $h$ the set $\Sigma$ has dimension $q-1$. For a complete proof 
see Proposition 2 in [Fo2]. The properness of $\pi|_\Sigma$ is easily 
satisfied by a small rotation of coordinates.) We may assume that $L=D\times B$ 
for some closed ball $B\subset \C^s$. To complete the proof we take 
$g=h\circ\psi$ where $\psi$ is a furnished by the following lemma.

\proclaim LEMMA 3.4: {\rm (Hypotheses as above)}
For every $\d>0$ there exists a \holo\ automorphism $\psi$ of $\C^n$ of 
the form $\psi(z,w)=(z,\b(z,w))$ such that $\psi(L)\cap \Sigma =\emptyset$
and $\sup_{x\in K} |\psi(x)-x| <\d$.

\demo Remark:
If $q=n$ then $\Sigma$ is a hypersurface in $\C^n$. In this case Lemma 3.4
fails since the complement $\C^n\bs \Sigma$ may be Kobayashi hyperbolic 
which would imply that any entire map $\psi\colon\C^n\to  \C^n\bs \Sigma$ 
has rank $<n$ at each point. 
\medskip

To prove Lemma 3.4 we shall need a version of Theorem 1.1 (or Theorem 2.1) 
from [FR] with a \holo\ dependence on parameters.  Recall that a vector field is 
{\it complete} if its flow exists for all times and all initial conditions. 
We shall consider holomorphic vector fields on $\C^n$ of the form
$$
	V(z,w)=\sum_{j=1}^s a_j(z,w){\di\over \di w_j} 
	\qquad (z\in\C^r,\ w\in\C^s),       			     \eqno(3.1)
$$
where the $a_j$'s are entire (or polynomial) functions on $\C^n=\C^r\times\C^s$.
Its flow remains in the level sets $\{z=const\}$, and $V$ is complete on $\C^n$ 
\iff\ $V(z,\cdotp)$ is complete on $\C^s$ for each $z\in \C^r$.

\proclaim LEMMA 3.5: If $s\ge 2$ then every polynomial vector field of type (3.1) 
on $\C^r\times \C^s$ is a finite sum of complete polynomial fields of the same type.

\demo Proof: 
We can write a polynomial field (3.1) as a finite sum $V(z,w)= \sum_{\a} z^\a V_\a(w)$ 
where  $V_{\a}(w)=\sum_{j=1}^s a_{\a,j}(w){\di\over \di w_j}$
and $z^\a=z_1^{\a_1}\cdots z_r^{\a_r}$. By [AL] every polynomial \hvf\ on $\C^s$ 
for $s\ge 2$ is a finite sum of complete polynomial fields (see the Appendix 
in [F1] for a short proof).  Hence each $V_{\a}(w)$ is a finite 
sum of complete polynomial fields. The products of such fields 
with $z^\a$ are complete on $\C^n$ which proves the result.

\medskip {\it Remark.}
For a more general result in this direction see Lemma 2.5 in [V] and
the recent preprint [Ku]. 
\medskip

Lemma 3.5 implies that the time-$t$ map of any entire \hvf\ (3.1) can be approximated, 
uniformly on any compact set on which it exists, by \holo\ automorphisms of 
$\C^n$ of the form $(z,w) \to \bigl(z,\varphi(z,w)\bigr)$ (Lemma 1.4 in [FR]). 
The same holds for time dependent entire \hvf s of the form (3.1). 
From this one obtains the following parametric version of Theorem 2.1 from [FR].

\proclaim COROLLARY 3.6:  
Assume that $\phi_t \colon \Omega_0\to\Omega_t$ $(t\in [0,T])$
is a smooth isotopy of biholomorphic maps between domains 
in $\C^n$, with $\phi_0$ the identity map on $\Omega_0$,
where $n=r+s$, $s\ge 2$, and each $\phi_t$ is of the form 
$\phi_t(z,w)=\bigl(z,\varphi_t(z,w) \bigr)$ $(z\in\C^r,\ w\in \C^s)$.
If $M\subset \Omega_0$ is a compact \pc\ set such that $\phi_t(M)$ is \pc\ 
in $\C^n$ for every $t\in [0,T]$ then $\phi_{T}$ can be approximated,
uniformly on $M$, by automorphisms of $\C^n$ of the form 
$(z,w)\to \bigl(z,\varphi(z,w)\bigr)$.

\demo Proof of Lemma 3.4:
Let $\Sigma$, $K$ and $L$ be as in the lemma. 
The set $\Sigma' :=\Sigma\cap (D\times\C^s)$ is \pc. 
Since $\Sigma\cap K=\emptyset$, $M :=K\cup \Sigma'$ is also \pc. 
Let $\theta_t(z,w) = (z,e^t w)$. Since the fibers $K_z$ ($z\in D$) 
are convex and contain the origin, the subvariety $\theta_t(\Sigma) \subset \C^n$ 
is disjoint from $K$ for every $t\ge 0$ and hence $M_t=K\cup \theta_t(\Sigma')$ 
is \pc. Clearly $\theta_T(\Sigma) \cap L=\emptyset$  for 
a sufficiently large $T>0$.

Consider the flow $\phi_t$ which rests on a \nbd\ of 
$K$ and equals $\theta_t$ on a \nbd\ of $\Sigma'$. Since $\phi_t(M)=M_t$ 
is \pc\ for every $t\ge 0$, Corollary 3.6 gives an automorphism 
$\phi(z,w) =\bigl(z,\varphi(z,w)\bigr)$ approximating 
$\phi_T$ uniformly on a \nbd\ of $M$. Thus $\phi$ is close to
the identity on a \nbd\ of $K$ and $\phi(\Sigma')\cap L=\emptyset$. 
Since $\phi$ maps each of the affine planes $\{z\}\times \C^s$
to itself, it follows that $\phi(\Sigma)\cap L=\emptyset$. 
The inverse $\psi=\phi^{-1}$ clearly satisfies Lemma 3.4.

\beginsection 4. Compositional splitting of biholomorphic mappings

Let $X$ be a complex manifold of dimension $n$. An injective 
\holo\ map $\gamma\colon V\to X$ in an open set $V\subset X$ will 
be called {\it biholomorphic}. Set $\disc=\{\z\in\C\colon |\z|<1\}$.
Suppose that $\cF$ is a nonsingular holomorphic foliation of $X$ of 
dimension $p$ and codimension $q=n-p$. Every $x\in X$ is contained
in a {\it distinguished chart} $(U,\phi)$, where $U\subset X$ 
is an open subset containing $x$ and $\phi\colon U\to \disc^n \subset\C^n$ 
is a biholomorphic map onto the open unit polydisc in $\C^n$ 
such that, in the coordinates $(z,w)$ on $\disc^n=\disc^p\times\disc^q$ 
($z\in \disc^p$, $w\in \disc^q$), $\phi(\cF|_U)$ is given by 
$\{w=c\}$, $c\in\disc^q$. Fix a number $0<r<1$. For any distinguished 
chart $(U,\phi)$ on $X$ let $U'\subset U$ be defined by 
$\phi(U')=(r\disc^p)\times \disc^q$.
Given any relatively compact set $V\ss X$, there exists 
a finite collection of distinguished charts
$\cU=\{(U_j,\phi_j)\colon 1\le j\le N\}$ 
such that $\bar V\subset \cup_{j=1}^N U'_j$ and 
$\cU$ is $\cF$-{\it regular} in the sense of Definition 1.5 
in [God, p.\ 72] (this means that for every 
$U_i,U_j\in\cU$ the set $\bar U_i\cap \bar U_j$ 
is contained in a distinguished chart).

%
%
\medskip DEFINITION.
A biholomorphic map $\gamma\colon V\to X$ 
is said to be an {\it $\cF$-map} if there exists $\cU$ as above
such that for every $(U_j,\phi_j)\in \cU$ the restriction of $\gamma$ to 
$V\cap U'_j$ has range in $U_j$ and is of the form $(z,w)\to (c_j(z,w),w)$ 
in the distinguished holomorphic coordinates on $U_j$.
\medskip

Thus an $\cF$-map preserves the leaves of $\cF$ and does not
permute the connected components of a global leaf intersected with 
any of the distinguished sets $U_j$. The definition is good since 
the transition map between a pair of distinguished charts 
preserves this form of the map. Any $\gamma$ preserving
the leaves of $\cF$ (in the sense that $x$ and $\gamma(x)$ belong
to the same leaf) which is close to the identity map in 
the {\it fine topology} on $X$, defined by $\cF$, 
is of this form. (The restriction of the fine topology to any 
distinguished local chart $U\simeq \disc^p\times \disc^q$ 
is the product of the usual topology on $\disc^p$ and the 
discrete topology on $\disc^q$. 
For further details see [God], pp.\ 2-3 and pp.\ 71--75.)

\proclaim THEOREM 4.1:
Let $A$ and $B$ be compact sets in a complex manifold $X$ such that
$D=A\cup B$ has a basis of Stein \nbd s in $X$ and 
$\bar{A\bs B}\cap \bar {B\bs A} =\emptyset$. Given an open set $\wt C \subset X$ 
containing $C:=A\cap B$ there exist open sets $A'\supset A$, $B'\supset B$, 
$C'\supset C$, with $C'\subset A'\cap B'\subset \wt C$, satisfying the following.
For every biholomorphic map $\gamma\colon \wt C\to X$ which is 
sufficiently uniformly close to the identity on $\wt C$ there exist 
biholomorphic maps $\a \colon A'\to X$, $\b\colon B'\to X$, 
uniformly close to the identity on their respective domains and satisfying 
$$
	\gamma = \beta\circ \alpha^{-1}	   \quad {\rm on\ } C'.
$$ 
If $\cF$ is a holomorphic foliation of $X$ and $\g$ is an $\cF$-map on $\wt C$ 
then we can choose $\a$, $\b$ to be $\cF$-maps on $A'$ resp.\ on $B'$. 
If $X_0$ is a closed complex subvariety of $X$ such that $X_0 \cap C=\emptyset$ 
then we can choose $\a$ and $\b$ as above such that they are tangent to the 
identity map to any given finite order along $X_0$.

Theorem 4.1, which is a key ingredient in our construction of noncritical
holomorphic functions and submersions, is proved in this section 
by a Kolmogorov-Nash-Moser type rapidly convergent iteration. 
(We shall only need it for the trivial foliation with $X$ as the only leaf,
but we prove the extended version for possible future applications.) It will be used in 
Propositions 5.2 and 6.1 below to patch pairs of noncritical functions or 
submersions (we have already explained in the introduction why the standard 
$\dibar$-theory does not suffices). Similar decompositions have been used 
in the theory of quasiconformal mappings and in complex dynamics. 
For example, a theorem of Pfluger [Pf] from 1961 asserts that 
every orientation preserving quasiconformal homeomorphism $\gamma\colon \R\to\R$ is 
the restriction to $\R$ of the composition $\beta \circ \alpha^{-1}$, where 
$\a, \b$ are conformal maps of the upper (resp.\ the lower) half plane 
to itself which map $\R$ to $\R$. (See also [LV, p.\ 92].)

We begin with preparatory  results. We fix once and for all a complete 
distance function $d\colon X\times X\to \R_+$ induced by a smooth Riemannian metric
on $TX$. Given a subset $A\subset X$ and an $r>0$ we set 
$$ 
	A(r) = \{x\in X\colon d(x,y)< r\ {\rm for\ some\ } y\in A\}.
$$
If $A$ is a (relatively) compact, smoothly bounded domain in $X$
then for all sufficiently small $r>0$ the set $A(r)$ is a smoothly 
bounded open domain. 

We say that the subsets $A,B\subset X$ are {\it separated\/} 
if $\bar{A\bs B}\cap \bar {B\bs A} =\emptyset$.

\proclaim LEMMA 4.2: Given $A,B\subset X$ and $r>0$ we have 
$(A\cup B)(r) = A(r) \cup B(r)$ and $(A\cap B)(r)\subset A(r)\cap B(r)$. 
If $A$ and $B$ are (relatively) compact and separated in $X$ 
then for all sufficiently small $r>0$ we also have $(A\cap B)(r)=A(r)\cap B(r)$ 
and the sets $\bar{A(r)}$, $\bar {B(r)}$ are separated.

\demo Proof: 
The first two properties are immediate. 
Now write $A=(A\bs B)\cup (A\cap B)$, $B=(B\bs A)\cup (A\cap B)$ 
and apply the first property to get
$$ 
	A(r)=(A\bs B)(r) \cup (A\cap B)(r),\quad B(r)=(B\bs A)(r) \cup (A\cap B)(r).
$$
If $\bar{A\bs B}\cap \bar {B\bs A} =\emptyset$ then for all sufficiently small 
$r>0$ we have $(A\bs B)(r) \cap (B\bs A)(r)=\emptyset$ (in fact, even the 
closures of $(A\bs B)(r)$ and $(B\bs A)(r)$ are disjoint). Hence the previous display 
gives $A(r)\cap B(r)= (A\cap B)(r)$ as well as the separation property for 
the pair $\bar{A(r)}$, $\bar{B(r)}$.

\proclaim  LEMMA 4.3: Let $A,B\subset X$ be compact sets in a complex manifold $X$
satisfying $\bar{A\bs B}\cap \bar {B\bs A} =\emptyset$. Assume that $A\cup B$ has a 
basis of Stein \nbd s. Given open sets $\wt A\supset A$, $\wt B\supset B$, 
$\wt C\supset C=A\cap B$, $\wt D\supset A\cup B$, there exist compact sets 
$A',B'\subset X$ satisfying the following:
\item{(a)} $A\subset A'\subset \wt A$, $B\subset B'\subset \wt B$, $A'\cap B'\subset \wt C$,
\item{(b)} $\bar{A'\bs B'}\cap \bar {B'\bs A'} =\emptyset$, 
\item{(c)} the set $D'=A'\cup B' \subset \wt D$ is the closure of a smoothly 
bounded \spsc\ Stein domain in $X$.

\demo Proof: If $r>0$ is chosen sufficiently small then by Lemma 4.2 we have
$A(r)\ss \wt A$, $B(r) \ss \wt B$, $A(r)\cap B(r) = C(r) \ss \wt C$,
and the sets $\bar{A(r)}$, $\bar{B(r)}$ are separated.  By assumption 
there is a closed \spsc\ Stein domain $D' \subset X$ with 
$A\cup B \subset D'\subset A(r)\cup B(r)$. It is easily verified 
that the sets $A'=\bar {A(r)} \cap D'$, $B'=\bar {B(r)}\cap D'$ 
satisfy the stated properties.
\medskip

Due to Lemma 4.3 it suffices to prove Theorem 4.1 under the 
assumption that $X$ is a Stein manifold, $A,B\subset X$ is  pair
of separated compact subsets and $D=A\cup B$ is the closure of a smoothly 
bounded strongly pseudoconvex domain. We assume this to be the case 
for the rest of this section.

Let $\cF$ be a  \holo\ foliation of $X$ with leaves $\cF_x$
$(x\in X)$. By Cartan's Theorem A the tangent bundle $T\cF \subset TX$ 
of $\cF$ is spanned by finitely many \hvf s $L_1,L_2,\ldots,L_m$ on $X$.
(We may have to shrink $X$ a bit.)
Denote by $\theta^j_t(x)$ the flow of $L_j$ for time $t\in \C$,
solving ${\di\over \di t}\theta^j_t(x)= L_j(\theta^j_t(x))$
and $\theta^j_0(x)=x$. The map $\theta^j$ is defined and \holo\ 
for $(x,t)$ in an open \nbd\ of $X\times \{0\}$ in $X\times \C$.
Their composition
$$
	\theta(x,t)=\theta(x,t_1,\ldots,t_m) :=
        \theta^m_{t_m}\circ\cdots\circ\theta^2_{t_2}\circ\theta^1_{t_1}(x) \in X 
$$
is a \holo\ map on an open \nbd\ $U\subset X \times\C^m$ of the zero section 
$X\times \{0\}^m$, satisfying $\theta(x,t)\in \cF_x$ for all $(x,t)\in U$ and 
$$
     \theta(x,0)=x,\qquad {\di\over \di t_j} \,\theta(x,t)|_{t=0}  = L_j(x) 
     \quad (x\in X,\ 1\le j\le m). 
$$
Hence $\Theta:= \di_t \theta|_{t=0}$ maps the trivial bundle $X\times\C^m$ 
surjectively onto the tangent bundle $T\cF$ of $\cF$. Splitting 
$X\times \C^m=E\oplus \ker \Theta$ we see that $\Theta \colon E \to T\cF$ 
is an isomorphism of holomorphic vector bundles. In any \holo\ vector bundle 
chart on $E$ we have a Taylor expansion
$$ 
	\theta(x,t_1,\ldots,t_m)= x+ \sum_{j=1}^m t_j\, L_j(x) + O(|t|^2)  \eqno(4.1)
$$
where the remainder $O(|t|^2)$ is uniform on any compact subset of the base set. 

Choose a Hermitean metric $|\cdotp|_E$ on $E$. Given an open set 
$V\subset X$ and a section $c \colon V\to E|_V$ we shall write 
$||c||_V=\sup_{x\in V}|c(x)|_E$.  By the construction of $\theta$ and 
$E$, $x\to \theta(x,c(x))$ is an $\cF$-map provided $||c||_V$  
is sufficiently small. 

Given a map $\gamma\colon V \to X$ we define 
$||\gamma-id||_V = \sup_{x\in V} d(\gamma(x),x)$, and we 
say that $\g$ is {\it $\e$-close to the identity on $V$} 
if $||\g-id||_V <\e$. The following lemma follows from the 
implicit function theorem.

\proclaim LEMMA 4.4: For every open relatively compact set $V\ss X$ there exist 
constants $M_1\ge 1$ and $\e_0>0$ satisfying the following property. For any 
$\cF$-map $\gamma\colon V\to X$ with $||\gamma-id||_V < \e_0$  
there is a unique \holo\ section $c\colon V\to E$ of $E|_V\to V$ such that 
for every $x\in V$ we have $\theta(x,c(x))= \gamma(x)$ and
$$
	M_1^{-1} |c(x)| \le d(\gamma(x),x) \le M_1 |c(x)|.
$$

If $\cF$ is the trivial foliation with $X$ as the only leaf, Lemma 4.4
asserts that every biholomorphic map $\gamma\colon V\to X$ sufficiently close 
to the identity map has the form 
$\gamma(x)=\theta(x,c(x))$ for some \holo\ section $c\colon V\to TV$.

We shall write the composition $\g\circ\a$ simply as $\g\a$. 
From now on all our sets in $X$ will be assumed contained in a fixed
relatively compact set for which Lemma 4.4 holds with a constant 
$M_1$. Recall that $V(\d)$ denotes the open $\d$-\nbd\ of $V\subset X$
with respect to the distance function $d$.

\proclaim LEMMA 4.5: 
Let $V\ss X$. There are constants $\d_0>0$ (small) and $M_2>0$ (large)
with the following property. Let $0< \d<\d_0$ and $0< 4\e <\d$. Assume that 
$\a,\b,\g \colon V(\d)\to X$ are $\cF$-maps which are $\e$-close 
to the identity on $V(\d)$. Then $\wt \g:= \b^{-1} \g \a \colon V\to X$ 
is a well defined $\cF$-map on $V$. Write 
$$
	\eqalign{ \a(x)=\theta(x,a(x)), & \quad \b(x)=\theta(x,b(x)), \cr
            \g(x)=\theta(x,c(x)),   & \quad
	\wt\g(x) =\theta(x, \wt c(x)), \cr}
$$
where $a,b,c$ are sections of $E|_{V(\d)} \to V(\d)$ and 
$\wt c$ is a section of $E|_V \to V$ given by Lemma 4.4. 
Then 
$$  
	||\,\wt c -(c+a-b)||_{V} \le  M_2\d^{-1} \e^2.         \eqno(4.2)   
$$
If $c=b-a$ on $V(\d)$ then $||\wt c||_V\le M_2\d^{-1} \e^2$
and $||\wt \g -id||_V \le M_1M_2\d^{-1} \e^2$.

\demo Proof:
The conditions imply that $\g\a$  maps $V$ biholomorphically
onto a subset of $V(2\e)$. Since $\b$ is $\e$-close to the identity map on $V(\d)$, 
the degree theory shows that its range contains $V(\d-\e)$. Hence
$\b^{-1}$ is defined on $V(\d-\e)$ and is $\e$-close to the 
identity on this set. Since $4\e<\d$, it follows that $\wt \g= \b^{-1}\g\a$ 
is defined on $V$ and maps $V$ biholomorphically onto 
a subset of $V(3\e)\ss V(\d)$. 

To prove the estimate (4.2) we choose a \holo\ vector bundle chart 
on $\pi\colon E\to X$ over an open set $\wt U\subset X$ and let $U\ss \wt U$.
We shall use the expansion (4.1) for $\theta$ on $\pi^{-1}(U) \subset E$; 
this suffices since $\bar {V(\d)}$ can be covered by finitely many such sets $U$.  
We replace the fiber variable $t$ in (4.1) by one of the functions $a(x),b(x)$, or $c(x)$. 
These are bounded on $V(\d)$ by a $M_1\e$ where $M_1$ is the constant from Lemma 4.4. 
This gives for $x\in U\cap V(\d)$:
$$ 
	\eqalign{ \a(x) &= x+\sum_{j=1}^m a_j(x)L_j(x)+ O(\e^2), \cr
              \b(x) &= x+\sum_{j=1}^m b_j(x)L_j(x)+ O(\e^2), \cr
              \g(x) &= x+\sum_{j=1}^m c_j(x)L_j(x)+ O(\e^2),\cr}
$$
where the remainder term $O(\e^2)$ is uniform with respect to 
$x\in U\cap V(\d)$. For $x\in U\cap V$ this gives 
$$ \eqalign{
  \g(\a(x)) &= \a(x) + \sum_{j=1}^m c_j(\a(x)) L_j(\a(x))+ O(\e^2) \cr
            &=     x + \sum_{j=1}^m \left( a_j(x)+c_j(x)\right) L_j(x)  \cr
            &\qquad + \sum_{j=1}^m \bigl( c_j(\a(x))L_j(\a(x)) -
            c_j(x)L_j(x) \bigr) + O(\e^2). \cr}
$$
To estimate terms in the last sum we fix $j$ and write  
$g(x)=c_j(x) L_j(x)$ for $x\in U\cap V(\d)$. Since $||c_j||_{V(\d)} < M_1\e$
and $4\e<\d$, the Cauchy estimates imply $||dc_j||_{U\cap V(\e)} = O(\e/\d)$ 
(here $dc_j$ denotes the differential of $c_j$). 
Since $L_j$ is \holo\ in a \nbd\ of $\bar{V(\d)}$, we may assume 
that its expression in the local coordinates on $U$ is uniformly bounded 
and has uniformly bounded differential. 
This gives $||dg||_{U\cap V(\e)} =O(\e/\d)$. 
Since $d(x,\a(x))<\e$, there is a smooth arc $\l\colon [0,1] \to U$, of 
length comparable to $\e$, such that $\l(0)=x$ and $\l(1)=\a(x)$. Then
$$ 
	|g(\a(x))-g(x)|\le \int_0^1 |dg(\l(\tau)| \cdotp |\l'(\tau)|\, d\tau
                  \le O(\d^{-1} \e^2)
$$
(the extra $\e$ is contributed by the length of $\l$). 
This gives for $x\in U\cap V$
$$ 
   \g(\a(x)) =  
   x+ \sum_{j=1}^m \left( a_j(x)+c_j(x)\right) L_j(x)  + O(\d^{-1} \e^2).
$$
The same argument holds for the composition of several maps
provided that $\e$ is sufficiently small in comparison to $\d$;
the error term remains of order $O(\d^{-1} \e^2)$. 

It remains to find the Taylor expansion of $\b^{-1}$ on the 
set $U\cap V(2\e)$ where $U$ is a local chart as above.
Set $\wt\b(x)= x-\sum_{j=1}^m b_j(x)L_j(x)$ for $x\in U\cap V(\d)$. 
Assuming that $\b(x)\in U\cap V(2\e)$ we obtain 
$$ 
	\eqalign{ \wt\b(\b(x)) &= \b(x)-\sum_{j=1}^m b_j(\b(x) L_j(\b(x)) \cr
                           &= \quad x+\sum_{j=1}^m 
      \bigl( b_j(x)L_j(x) - b_j(\b(x)) L_j(\b(x)) \bigr) +O(\e^2) \cr
                           &= \quad x+O(\d^{-1} \e^2). }
$$
We have estimated the terms in the parentheses on the middle line by 
$O(\d^{-1} \e^2)$ in exactly the same way as above, using the Cauchy estimates
and integrating over an arc of length comparable to $\e$.
Writing $\b(x)=y\in U\cap V(2\e)$, $x=\b^{-1}(y)$, the above gives 
$\wt \b(y)=\b^{-1}(y)+O(\d^{-1} \e^2)$ and therefore
$$  
	\b^{-1}(y)= y-\sum_{j=1}^m b_j(y)L_j(y) + O(\d^{-1} \e^2). 
$$
The same argument as before gives 
$$ 
	\wt \g(x)= (\b^{-1}\gamma\alpha)(x) 
	         = x+\sum_{j=1}^m 
	\bigl( c_j(x)+a_j(x)-b_j(x) \bigr) L_j(x) + O(\d^{-1} \e^2)
$$
for $x\in U\cap V$. This proves the estimate (4.2). 
\medskip

{\it Remark.} The proof of Lemmas 4.4 and 4.5 shows that for each fixed open set
$V_0\ss X$ the constants $M_1,M_2,\d_0$ may be chosen independent of $V$ 
for any open set $V\subset V_0$. In this case $O(\d^{-1}\e^2)$ means
$\le C\d^{-1}\e^2$, with $C$ independent of $\e$, $\d$ and $V$.

\proclaim 4.6 LEMMA:   
Let $E\to X$ be a holomorphic vector bundle over a Stein manifold $X$. 
Let $U, V \subset X$ be open sets such that $\bar{U\bs V} \cap \bar{V\bs U} =\emptyset$
and $D=U\cup V$ is a relatively compact, smoothly bounded, \spsc\ domain in $X$.
Set $W=U\cap V$. There is a constant $M_3 \ge 1$ such that for every bounded 
\holo\ section $c\colon W\to E|_W$ there exist bounded \holo\ sections 
$a\colon U\to E|_U$, $b\colon V\to E|_V$ satisfying 
$$ 
	c=b|_W -a|_{W},\quad ||a||_U< M_3 ||c||_W,\ \ ||b||_V< M_3 ||c||_W. 
$$   
Such $a$ and $b$ are given by bounded linear operators between the spaces 
of bounded \holo\ sections of $E$ on the respective sets. The constant $M_3$ can be 
chosen uniform for all such pairs $(U,V)$ in $X$ close to an initial pair $(U_0,V_0)$ 
provided that $D=U\cup V$ is sufficiently $\cC^2$-close to $D_0=U_0\cup V_0$.
If $X_0$ is a closed complex subvariety of $X$ and $X_0\cap \bar W=\emptyset$
then for every $s\in\N$ we can insure in addition that $a$ and $b$ vanish 
to order $s$ on $X_0$.

\demo Proof: This is a standard application of the solvability of the 
$\dibar$-equation. We give a brief sketch for the sake of completeness.

Condition (b) insures that there is a smooth function $\chi\colon X\to [0,1]$ 
which equals zero in a \nbd\ of $\bar {U\bs V}$ and equals one in a \nbd\ of 
$\bar {V\bs U}$. Since $D=U\cup V$ is a relatively compact \spsc\ domain in $X$,
there exists a bounded linear solution operator $T$ for the $\dibar$-equation
associated to sections of $E\to X$ over $D$. Precisely, for any bounded
$\dibar$-closed $E$-valued $(0,1)$-form $g$ on $D$ we have 
$\dibar_E(Tg)=g$ and $||Tg||_D \le {\rm const}||g||_D$ and the  
constant can be chosen uniform for all domains in $X$ which are sufficiently 
$\cC^2$-close to an initial \spsc\ domain. (For functions this can be found 
in [HL1, p.\ 82]; the problem for sections of a vector bundle $E$ can be reduced 
to that for functions by embedding $E$ as a subbundle of a trivial bundle over $X$.)

Observe that $\chi c$ extends to a bounded smooth section of $E$ over $U$
and $(\chi-1)c$ extends to a bounded section over $V$.
Since $\supp (\dibar \chi)\cap D \subset W=U\cap V$, the bounded
$(0,1)$-form $g=\dibar(\chi c)=\dibar((\chi-1)c)= c \dibar \chi$
on $W$ extends to a bounded $(0,1)$-form on $D$ which is zero
outside of $W$. It is immediate that the pair of sections
$$ a = -\chi c    + T(g)|_U, \qquad
   b = (1- \chi)c + T(g)|_{V}
$$
satisfies Lemma 4.6. The last statement (regarding the 
interpolation on $X_0$) follows in the case of 
functions from [FP3, Lemma 3.2]; the same proof applies to sections 
of $E\to X$ by embedding $E$ into a trivial bundle over $X$.

\proclaim LEMMA 4.7:
Let $A,B\subset X$ be compact sets such that 
$\bar{A\bs B} \cap \bar{B\bs A} =\emptyset$ and  
$D=A\cup B$ is a closed, smoothly bounded, \spsc\ domain in $X$. 
Let $\cF$ be a  \holo\ foliation of $X$ 
and let $X_0$ be a closed complex subvariety of $X$
with $X_0\cap C=\emptyset$, where $C=A\cap B$. 
Then there are constants $r_0>0,\d_0>0$ (small)
and $M_4,M_5 > 1$ (large) satisfying the following. 
Let $0<r \le r_0$, $0<\d \le \d_0$ and $s\in \N$. For every $\cF$-map 
$\g\colon C(r+\d)\to X$ satisfying 
$4M_4 ||\g-id||_{C(r+\d)}< \d$ there exist $\cF$-maps 
$\a\colon A(r+\d)\to X$ and $\b\colon B(r+\d)\to X$, tangent to the 
identity map to order $s$ along $X_0$, such that $\wt\g :=\b^{-1}\g\a$ 
is an $\cF$-map on $C(r)$ satisfying 
$$
	||\wt\g -id||_{C(r)}  
	< M_5 \d^{-1} ||\g-id||_{C(r+\d)}^2.      \eqno(4.3)
$$

\ni\it Proof. \rm
If $r_0$ and $\d_0$ are chosen sufficiently small,
the set $D(t)$ is a small $\cC^2$-perturbation of the \spsc\ domain 
$D=A\cup B$ for every $t\in [0,r_0+\d_0]$ and hence we can use the same 
constant as a bound on the sup-norm of an operator solving the 
$\dibar$-problem on $D(t)$. 

Let $\e=||\g-id||_{C(r+\d)}$. By Lemma 4.4 there is a 
\holo\ section $c\colon C(r+\d)\to E$, 
with $||c||_{C(r+\d)}\le M_1\e$, such that $\g(x)=\theta(x,c(x))$. 
(Here we can use the constant $M_1$ for the set $D(r_0+\d_0)$.)
Write $c=a-b$ where $a$ is a section of $E$ over 
$A(r+\d)$ and $b$ is a section of $E$ over $B(r+\d)$ furnished by 
Lemma 4.6. The sup-norms of $a$ and $b$ on their respective 
domains are bounded by $M_1M_3 \e$, where the constant $M_3$ from Lemma 4.5
can be chosen independent of $r$ and $\d$. Set
$$ 
	\eqalign{
	\a(x) &= \theta(x,a(x)) \qquad (x\in A(r+\d)), \cr
	\b(x) &= \theta(x,b(x)) \qquad (x\in B(r+\d)). \cr}
$$
By Lemma 4.4 we have $||\a-id||_{A(r+\d)} < M_1^2M_3\e$
and $||\b-id||_{B(r+\d)} < M_1^2M_3\e$. Set $M_4=M_1^2M_3$.
If $0<4M_4\e<\d$ then by Lemma 4.5 the composition $\wt \g=\b^{-1}\g\a$ is 
an $\cF$-map on $C(r)$ satisfying the estimate (4.3) with 
$M_5= M_2M_4^2=M_1^4 M_2 M_3^2$.
This completes the proof.

\demo  Proof of Theorem 4.1:
By Lemma 4.3 we may assume that $D=A\cup B$ is the closure of a smooth 
strongly \psc\ domain in $X$ and $\bar{A\bs B}\cap \bar{B\bs A}=\emptyset$.
Choose a sufficiently small number $0< r_0 <1$ such that the initial 
$\cF$-map $\g$ is defined on the set $C_0 := C(r_0)$ and Lemma 4.7 holds for 
all $\d,r>0$ with $\d+r\le r_0$. For each $k=0,1,2,\ldots$ we set  
$$
	r_k=r_0\prod_{j=1}^k (1-2^{-j}),\quad 
        \d_k = r_k - r_{k+1} = r_k 2^{-k-1}.  
$$
The sequence $r_k>0$ is decreasing, $r^* = \lim_{k\to \infty}r_k > 0$,
$\d_k > r^* 2^{-k-1}$ for all $k$, and $\sum_{k=0}^\infty \d_k= r_0-r^*$. 
Set $A_k=A(r_k)$, $B_k=B(r_k)$, $C_k=C(r_k)$. 
We choose $r_0>0$ sufficiently small such that $C_k=A_k\cap B_k$ 
for all $k$ (Lemma 4.2).  

Let $\e_0 :=||\gamma -id||_{C_0}$. Assuming that 
$4M_4 \e_0 <\d_0=r_0/2$, Lemma 4.7 gives $\cF$-maps 
$\a_0\colon A_0\to X$ and $\b_0\colon B_0\to X$ such that
$\g_1=\b_0^{-1}\g\a_0 \colon C_1\to X$ is an $\cF$-map 
defined on $C_1$, satisfying 
$$ 
	||\g_1-id||_{C_1} < M_5\d_0^{-1} \e^2_0 < 2M\e_0^2,
$$
where we have set $M=M_5/r^*$. 
Define $\e_1=||\g_1-id||_{C_1}$, so $\e_1<2M\e_0^2$. 
Assuming for a moment that $4M_4\e_1 <\d_1$, we can apply
Lemma 4.7 to obtain a pair of $\cF$-maps $\a_1\colon A_1\to X$,
$\b_1\colon B_1\to X$ such that $\g_2=\b_1^{-1}\g_1 \a_1\colon C_2\to X$
is an $\cF$-map satisfying 
$$
	\e_2:= ||\g_2-id||_{C_2} < M_5\d_1^{-1} \e^2_1 < 2^2 M \e_1^2.
$$
Continuing inductively we obtain sequences of $\cF$-maps
$$
	\a_k\colon A_k\to X,\quad \b_k \colon B_k\to X, 
	\quad \g_k\colon C_k\to X
$$
such that $\g_{k+1}=\b_k^{-1}\g_k\a_k \colon C_{k+1}\to X$ is an 
$\cF$-map satisfying 
$$ 
	\e_{k+1}:= ||\g_{k+1}-id||_{C_{k+1}} 
	< M_5\d_k^{-1} \e^2_k < 2^{k+1} M \e_k^2.       \eqno(4.4)
$$
The necessary condition for the induction step is that
$4M_4\e_k<\d_k$ holds for each $k$. Since $\d_k > r^* 2^{-k-1}$, 
it suffices to have  
$$ 
	4M_4 \e_k< r^* 2^{-k-1}  \qquad (k=0,1,2,\ldots). \eqno(4.5)
$$
In order to obtain convergence of this process we need the following.

\proclaim LEMMA 4.8: Let $M,M_4\ge 1$. Let the sequence $\r_k>0$ be defined 
recursively by $\r_0=\e_0>0$ and $\r_{k+1}= 2^{k+1} M \r_k^2$ for $k=0,1,\ldots$.
If $\e_0< r^*/32 MM_4$ then $\r_k < (4M\e_0)^{2^k} < (1/8)^{2^k}$ and
$4M_4\r_k < r^* 2^{-k-1}$ for all $k=0,1,2,\ldots$.

Assuming Lemma 4.8 we complete the proof
of Theorem 4.1 as follows. From (4.4) we see that $\e_k  \le \r_k$
where $\r_k$ is the sequence from Lemma 4.8. From the 
assumption $\e_0< r^*/32 MM_4$ 
we obtain $q:=4M\e_0 < r^*/8M_4 < 1/8$ (since $0<r^*<1$ and $M_4\ge 1$). 
Hence the sequence $\e_k=||\g_k-id||_{C_k} < q^{2^k} < (1/8)^{2^k}$ 
converges to zero very rapidly as $k\to\infty$.
The second estimate on $\r_k$ in Lemma 4.8 insures that (4.5) holds
and hence the induction described above works.  

Setting $\wt\a_k=\a_0\a_1\cdots\a_k \colon A_k\to X$,
$\wt \b_k=\b_0\b_1\cdots\b_k \colon B_k\to X$, we have 
$\g_{k+1} = \wt \b_k^{-1} \g \wt \a_k$ on $C_{k+1}$
for $k=0,1,2,\ldots$. Our construction insures that,
as $k\to\infty$, the sequences $\wt\a_k$ resp.\ $\wt \b_k$ converge, 
uniformly on $A(r^*)$ resp.\ on $B(r^*)$, 
to $\cF$-maps $\a\colon A(r^*)\to X$ resp.\ $\b\colon B(r^*)\to X$. 
Furthermore, the sequence $\g_k$ converges uniformly 
on $C(r^*)$ to the identity map according to (4.4) 
and Lemma 4.8. In the limit we obtain $\b^{-1} \g \a =id$ on $C(r^*)$,
and hence $\gamma=\beta\alpha^{-1}$ on $\a(C(r^*))$.
If $\e_0>0$ is chosen sufficiently small (for a fixed $r_0$)
then the latter set contains a \nbd\ $C'$ of $C$.
This completes the proof of Theorem 4.1, provided that Lemma 4.8 holds.

\demo Proof of Lemma 4.8:
The sequence is of the form $\r_k = 2^{a_k}M^{b_k}\e_0^{c_k}$
where the exponents satisfy the recursive relations
$$ \eqalign{
	a_{k+1} &= 2a_k+k+1,\quad a_0=0; \cr
	b_{k+1} &= 2b_k +1, \quad \quad \quad b_0=0; \cr
	c_{k+1} &= 2c_k, \quad \quad\quad\quad\ \   c_0=0.\cr}
$$
The solutions are $a_k=2^k\sum_{j=1}^k j2^{-j} < 2^{k+1}$,
$b_k=2^{k}-1$, $c_k=2^k$. Thus 
$$
	\r_k< 2^{2^{k+1}}M^{2^k}\e_0^{2^k} = 
	(4M\e_0)^{2^k}
$$
which proves the first required estimate. 
From the assumption $\e_0< r^*/32 MM_4$ we get $q:=4M\e_0 < r^*/8M_4<1/8$.
Hence $\rho_k < q^{2^k}< (1/8)^{2^k}$ and 
$$
	4M_4\r_k < (4M_4 q) q^{2^k-1} = 
	(4M_4\cdotp 4M\e_0) (1/8)^{2^k-1}< (r^*/2) 2^{-k} =r^* 2^{-k-1}
$$
for all $k\ge 0$ which proves the second estimate. Lemma 4.8 is proved.

\medskip
{\it Remark.}
The above construction actually gives nonlinear operators $\cA$, $\cB$ 
on the set of $\cF$-maps $\g$ which are sufficiently uniformly 
close to the identity on a fixed \nbd\ of $C$ such that the pair
of $\cF$-maps $\a=\cA(\g)$, $\b=\cB(\g)$ satisfies $\g=\b\a^{-1}$. 
This yields the analogous result for families 
$\{\gamma_p \colon p\in P\}$ of $\cF$-maps which depend continuously 
on a parameter $p$ in a compact Hausdorff space $P$ and which are 
sufficiently close to the identity map on a \nbd\ of $C$.

\beginsection 5. Construction of noncritical holomorphic functions

A compact set $K$ in a complex manifold $X$ is said 
to be a {\it  Stein compactum} if it has a basis of open Stein \nbd s.
Let $d$ be a distance function on $X$ induced by a smooth Riemannian metric on $TX$.
We shall use the terminology introduced in Section 4. Recall that
$||\g-id||_V=\sup_{x\in V} d(\g(x),x)$.

\proclaim LEMMA 5.1: Let $K$ be a Stein compactum in a complex manifold $X$.
Let $U\subset X$ be an open set containing $K$ and $f\colon U\to\C^q$
a \holo\ submersion for some $q\le \dim X$. Then there exist
constants $\e_0>0, M>0$ and an open set $V\subset X$, with
$K\subset V\subset U$, satisfying the following property. 
Given $\e\in (0, \e_0)$ and a  \holo\ submersion 
$g\colon U\to \C^q$ with $\sup_{x\in U} |f(x)-g(x)| < \e$ there
is a biholomorphic map $\g\colon V\to X$ satisfying 
$f=g\circ\g$ on $V$ and $||\g-id||_V < M\e$.

\demo Proof: We may assume that $U$ is Stein. Hence  
$TX|_U=\ker df \oplus E$ for some trivial rank $q$ holomorphic
subbundle $E\subset TX|_U$. Thus $E$ is spanned by 
$q$ independent \holo\ vector fields on $U$.
Denote by $\theta(x,t_1,\ldots,t_q)$ the composition of their local 
flows (see the construction of $\theta$ (4.1) in Sect.\ 4).
The map $\theta$ is defined in an open set $\Omega \subset U\times \C^q$
containing $U\times\{0\}^q$. For $x\in U$ write 
$\Omega_x=\{t\in\C^q \colon (x,t)\in \Omega\}$. After shrinking $\Omega$ 
we may assume that for each $x\in U$ the fiber $\Omega_x$ is connected and 
$F_x := \{\theta(x,t)\colon t\in \Omega_x\} \subset X$ is a local complex 
submanifold of $X$ which intersects the level set $\{f=f(x)\}$ transversely at $x$
(since $T_x F_x =E_x$ is complementary to the kernel of $df_x$).
By the implicit function theorem we may assume that (after shrinking $\Omega$)
the map $t\in \Omega_x\to f(\theta(x,t)) \in \C^q$ maps $\Omega_x$ 
biholomorphically onto a \nbd\ of the point $f(x)$ in $\C^q$. The same holds 
for the map $t\in\Omega_x \to g(\theta(x,t))$ provided that $g\colon U\to\C^q$ is 
sufficiently uniformly close to $f$ and we restrict $x$ to a compact subset 
of $U$. It follows that, if $V\ss U$ and $g$ is sufficiently
close to $f$ on $U$, there is for every $x\in V$ a unique 
point $c(x)\in \Omega_x$ such that $g(\theta(x,c(x)))=f(x)$.
Clearly $c\colon V\to \C^q$ is \holo\  and the map 
$\g(x)=\theta(x,c(x))\in X$ ($x\in V$) satisfies Lemma 5.1.

%
%
\medskip DEFINITION.
An ordered pair of compact sets $(A,B)$ in an $n$-dimensional
complex manifold $X$ is said to be a {\it special Cartan pair} if 
\item{(i)}  the sets $A$, $B$, $C:=A\cap B$, $A\cup B$ 
are Stein compacta (see above);  
\item{(ii)}  $\bar{A\bs B} \cap \bar{B\bs A} =\emptyset$, and 
\item{(iii)} there is an open set $U\supset B$ and an 
injective holomorphic map $\psi\colon U\to \C^n$ such that 
$\psi(C) \subset \C^n$ is \pc. 
\medskip

The following is the main step in the proof of Theorem 2.1.

\proclaim PROPOSITION 5.2: Let $(A,B)$ be a special Cartan pair in 
a complex manifold $X$ and let $f \in\cO(A)$ be a function whose critical 
set $P$ is finite and does not meet $C:=A\cap B$. Given $\e>0$ there 
exists $\wt f \in \cO(A\cup B)$ with the same critical set $P$ 
such that $\sup_A |\wt f-f|<\e$. If $X_0$ is a closed 
complex subvariety of $X$ with $X_0\cap C =\emptyset$ then for any 
$r\in \N$ we can choose $\wt f$ as above such that $\wt f-f$ 
vanishes to order $r$ on $X_0 \cap A$. In particular,
if $f$ is noncritical on $A$ then $\wt f$ is noncritical on $A\cup B$.

\demo Proof: 
We use the notation from (iii) in the definition of a special Cartan pair.
The function $f'=f\circ\psi^{-1}$ is defined and noncritical 
in an open set $\wt C\subset \C^n$ containing $\psi(C)$.  
Choose a compact \pc\ set $K$ with 
$\psi(C)\subset {\rm int} K \subset K \subset\wt C$.
By Theorem 3.1 we can approximate $f'$ uniformly on $K$ by a noncritical 
\holo\ function $g'\in\cO(\C^n)$. Thus $g=g'\circ\psi$ is noncritical
in a \nbd\ of $B$ and it approximates $f$ uniformly in a \nbd\ of $C$. 
If the approximation is sufficiently close then by Lemma 5.1 there is a 
biholomorphic map $\g$, uniformly close to the identity map in a \nbd\ 
of $C$, satisfying $f=g\circ \g$. By Theorem 4.1 we have 
$\g=\b\circ \a^{-1}$, where $\a$ is a biholomorphic map close to 
the identity in a \nbd\ of $A$ and $\b$ is a map with the analogous 
properties in a \nbd\ of $B$. Furhermore we insure that 
$\a$ agrees with the identity map to a sufficiently high order at each 
point of $(X_0\cup P)\cap A$. From $f=g\circ \g=g\circ\b\circ\a^{-1}$ 
(which holds in a \nbd\ of $C$) we obtain $f\circ \a=g\circ \b$. 
The two sides define a \holo\ function 
$\wt f \in\cO(A\cup B)$ which the stated properties.

%
%
\demo Proof of Theorem 2.1: 
We first consider the simplest case when $f$ is noncritical on $U$ and
$X_0=\emptyset$. By Corollary 2.8 in [HL3] there is a sequence of 
compact $\cO(X)$-convex subsets 
$A_0\subset A_1\subset \cdots\subset \cup_{k=0}^\infty A_k=X$
such that 

\smallskip
\item{(i)}  $K\subset {\rm Int A_0} \subset A_0 \subset U$, and
\item{(ii)} for every $k=0,1,2,\ldots$ we have $A_{k+1}=A_k\cup B_k$ 
where $(A_k,B_k)$ is a special Cartan pair in $X$.
\smallskip

Fix $\e>0$. Write $f_0=f$, $A_{-1}=K$.
Choose a sufficiently small number $\e_0\in (0, \e/2)$ 
such that every $g\in \cO(A_0)$ with 
$\sup_{A_0} |g-f_0|< 2\e_0$ is noncritical on $K$. 
Proposition 5.2 gives a noncritical function $f_1\in \cO(A_1)$
satisfying $\sup_{A_0} |f_1-f_0| < \e_0 < \e/2$. Now choose 
$\e_1 \in (0,\e_0/2)$ such that every function $g\in \cO(A_1)$ 
with $\sup_{A_1} |g-f_1|< 2\e_1$ is noncritical on $A_0$. 
Proposition 5.2 gives a noncritical function $f_2\in \cO(A_2)$ such 
that $\sup_{A_1} |f_2-f_1| < \e_1 < \e/4$. Continuing inductively 
we obtain a sequence of noncritical functions $f_k\in \cO(A_k)$ 
and a decreasing sequence $\e_k > 0$ with $\sum_{k=0}^\infty \e_k < \e$ such that 
$\sup_{A_k} |f_{k+1}-f_k| < \e_{k} < \e 2^{-k-1}$ for every
$k=0,1,2,\ldots$. The sequence $f_k$ converges uniformly on compacts 
in $X$ to $\wt f\in\cO(X)$ satisfying 
$\sup_K |\wt f -f|< \e$ and $\sup_{A_k} |\wt f-f_k|< 2\e_k$
for every $k=0,1,2,\ldots$. By the choice of $\e_k$ 
this insures that $\wt f$ is noncritical on $A_{k-1}$.
Since this holds for every $k$, $\wt f$ is noncritical on $X$.

Consider now the general case. Let 
$P=\{p_1,p_2,\ldots,\}$ denote the (discrete) critical set of $f\in \cO(U)$.
We replace $X_0$ by $X_0\cup P$. For each $j\in \N$ we choose a sufficiently 
large integer $n_j\in\N$ such that for every germ of a holomorphic function 
$g$ which vanishes to order $n_j$ at $p_j$ the germ of $f+g$ still 
has an isolated critical point at $p_j$.  In the sequel we shall often use the 
following elementary fact. Given a pair of compact sets $K\subset L$ in the domain of $f$,
with $K \subset {\rm Int}L$, we can choose $\eta>0$ such that for every 
$g\in\cO(L)$ which vanishes to order $n_j$ at every point $p_j\in P\cap K$ 
and satisfies $\sup_L |g|<\eta$ the critical set of $f+g$ in $K$ equals $P\cap K$. 

Denote by $\cJ\subset \cO_X$ the coherent analytic sheaf of ideals 
consisting of all germs of \holo\ functions on $X$ which vanish 
to order $r$ on $X_0$ and to order $n_j$ at $p_j\in P$ for 
every $j\in\N$. We can replace $f$ by a function holomorphic on $X$ 
such that the difference of the two functions
is a section of $\cJ$ near $X_0$ and is uniformly small 
on $K$ (see Lemma 8.1 in [FP1]). The new function (which we still 
denote $f$) may have additional critical points, but there is a 
\nbd\ $U\supset X_0\cup K$ such that ${\rm Crit}(f;U)=P$. 
Choose a compact $\cO(X)$-convex set $L\subset X$ containing $K$ 
in its interior. Fix an $\eta>0$. We claim that there exists 
an $f'\in\cO(L)$ satisfying the following:

\smallskip
\item{(i)}   ${\rm Crit}(f';L)=P\cap L$, 
\item{(ii)}  $f'-f$ is a section of $\cJ$ over $L$, and 
\item{(iii)} $|f'-f| <\eta$ on $K$.
\smallskip

Proof: By Lemma 8.4 in [FP2] there is a finite sequence 
$A_0\subset A_1\subset \cdots\subset A_{k_0}=L$ of compact $\cO(X)$-convex 
subsets such that for each $k=0,1,\ldots,k_0-1$ we have $A_{k+1}=A_k\cup B_k$, 
where $(A_k,B_k)$ is a special Cartan pair in $X$ and 

\smallskip
\item{(a)}  $K\cup (X_0\cap L) \subset A_0\subset \subset U$, 
\item{(b)}  $B_k\cap X_0=\emptyset$ for $k=0,1,\ldots,k_0-1$.
\smallskip

\ni (Our notation differs from [FP2]: the set $A_k$ in [FP2] is denoted 
$B_{k-1}$ in this paper, while the set $A_k$ in this paper is the same
as $\cup_{l=0}^k A_l$ in [FP2].) 
Assume inductively that for some $k<k_0$ we already have a function 
$f_k \in\cO(A_k)$ satisfying the above properties (i)--(iii) 
(with $f'$ replaced by $f_k$). Since $B_k \cap X_0 =\emptyset$, 
$f_k$ is noncritical in a \nbd\ of $A_k\cap B_k$ and hence Proposition 5.2 
furnishes a function $f_{k+1} \in\cO(A_{k+1})$ satisfying (i)--(iii) 
on its domain. After $k_0$ steps we obtain the desired function $f'\in\cO(L)$
thus proving the claim. 

\smallskip
In order to complete the induction step we show that there exists 
$h\in\cO(X)$ such that $h-f$ is a section of $\cJ$, $h$ approximates $f'$ 
uniformly on $L$, and there is a \nbd\ $\wt U\supset X_0\cup L$ such that 
${\rm Crit}(h;\wt U)=P$. By Cartan's Theorem A the sheaf $\cJ$ is finitely 
generated on the compact set $L$, say by functions $\xi_l\in\cO(X)$ $(l=1,2,\ldots,m)$. 
Since $f' - f$ is a section of $\cJ$ over a \nbd\ of $L$, we have 
$f' = f + \sum_{j=1}^m \xi_j g_j$ for some $g_j\in \cO(L)$. 
Since $L$ is $\cO(X)$-convex, we can approximate $g_j$ uniformly on 
a \nbd\ of $L$ by $\wt g_j\in\cO(X)$. The function 
$h=f + \sum_{j=1}^m \xi_j \wt g_j \in \cO(X)$ satisfies the stated properties 
provided that the approximation of $g_j$ by $\wt g_j$ was sufficiently 
close for every $j$. 

Note that $h$ satisfies the same properties on a \nbd\ of $L\cup X_0$
as $f$ did on a \nbd\ of $K\cup X_0$. The proof of Theorem 2.1 is completed 
by an obvious induction over a sequence of compact $\cO(X)$-convex sets 
$L_1\subset L_2\subset\cdots$ exhausting $X$ 
(compare with the noncritical case given above).

%
%
%
%
\beginsection 6. Construction of holomorphic submersions

In this section we prove Theorems 2.5 and 2.6 and Corollaries 2.10 and 2.11.
We begin with Theorems 2.5 and 2.6. Since the proof is fairly long,
we first explain the outline and then treat each of the 
main ingredients in a separate subsection.  

We are given a $q$-coframe $\theta=(\theta_1,\ldots,\theta_q)$ 
on $X$ such that $\theta|_U=df$ in an open set $U\supset K$ where 
$f\colon U\to\C^q$ is a holomorphic submersion. Our task is to 
deform $\theta$ to the differential $d\wt f$ where $\wt f\colon X\to \C^q$ 
is a holomorphic submersion which approximates $f$ uniformly on $K$. 
(We shall deal with interpolation along a subvariety $X_0\subset X$ 
in Subsect.\ 6.5.)

Let $\rho\colon X\to\R$ be a smooth \spsh\ Morse exhaustion function 
such that $\rho<0$ on $K$ and $\rho>0$ on $X\bs U$ ([H\"o2], Theorem 5.1.6.). 
Each sublevel set $\{\rho\le c\}$ is compact and $\cO(X)$-convex;
if $c\in \R$ is a regular value of $\rho$ then $\{\rho<c\}$ is a smooth
\spsc\ domain. The set of critical values of $\rho$ is discrete in $\R$ 
and hence at most countable, and each critical level contains a unique 
critical point. 

It suffices to explain how to approximate a submersion $f$ defined in a \nbd\ 
of $\{\rho \le  c_0\}$ (and with $df$ homotopic to $\theta$ through $q$-coframes) 
by a submersion $\wt f$ with similar properties defined in \nbd\ of $\{\rho\le c_1\}$, 
where $c_0< c_1$ is any pair of regular values of $\rho$. The construction is then 
completed by an obvious induction as in Theorem 2.1. Using a smooth cut-off
function in the parameter of the homotopy from $df$ to $\theta$ we can 
deform the $q$-coframe $\theta$ at each step to insure that 
$\theta=df$ in a \nbd\ of $\{\rho\le c_0\}$.

The construction of the extension breaks into two distinct arguments: 
(i) going through noncritical values of $\rho$ (mainly complex analysis), 
(ii) crossing a critical value (mainly topology and `convex analysis').

If $\rho$ has no critical values in $[c_0,c_1]$ then $\{\rho\le c_1\}$ is 
obtained from $\{\rho\le c_0\}$ by finitely many attachings of small 
{\it convex bumps}. In each step we approximately extend $f$ over the 
bump by Proposition 3.3, and we patch the two pieces using Theorem 4.1. 
In finitely many steps we obtain a submersion $f$ in a 
\nbd\ of $\{\rho\le c_1\}$ (Subsect.\ 1). 

Crossing a critical value of $\rho$ relies on a combination of three techniques:

\smallskip
\item{--}  smooth extension across a handle attached to $\{\rho\le c_0\}$
(Subsect.\ 2),
\item{--}  approximation by a \holo\ submersion defined in a \nbd\ of 
a handlebody (Subsect.\ 3), and 
\item{--} applying the noncritical case with a different \spsh\ function 
to extend across the critical level of $\rho$ (Subsect.\ 4).
\smallskip

The proof of Theorem 2.5 is completed immediately after Lemma 6.7 
in Subsect.\ 4, with the exception of the interpolation along a 
subvariety $X_0\subset X$ which is explained in Subsect.\ 5. 
There we also prove Corollaries 2.10 and 2.11.

%
%
%
%
\medskip
{\it 1. The noncritical case.}
A compact set $\wt A\subset X$ in a complex manifold $X$ is a 
{\it noncritical \spsc\ extension} of a compact set $A\subset \wt A$
if there is a smooth \spsh\ function $\rho$ in an open set 
$\Omega \supset \bar{\wt A\bs A}$ which has no critical points 
on $\Omega$ and satisfies
$$
	A\cap \Omega     = \{x\in \Omega \colon \rho(x)\le 0\},\quad   
        \wt A\cap \Omega = \{x\in \Omega \colon \rho(x)\le 1\}.
$$
Note that for each $t\in [0,1]$ the set 
$A_t=A\cup \{\rho\le t\} \subset X$ is a smooth (closed) \spsc\ domain 
in $X$, and the family smoothly increases from $A=A_0$ to $\wt A=A_1$. 
We say that {\it $X$ is a noncritical \spsc\ extension of $A$} if 
there exists a smooth exhaustion function $\rho\colon X \to \R$
such that $A=\{\rho\le 0\}$ and $\rho$ is \spsh\ and without
critical points on $\{\rho\ge 0\}=X\bs {\rm int A}$. 

\proclaim PROPOSITION 6.1:
Let $X$ be a Stein manifold and $\wt A\subset X$ a noncritical \spsc\ 
extension of $A\subset \wt A$. If $f\colon A\to \C^q$ is a \holo\ 
submersion with $q<\dim X$ then for every $\e>0$ there exists a 
\holo\ submersion $\wt f\colon \wt A \to \C^q$ satisfying 
$\sup_A |f-\wt f|<\e$. 

\proclaim COROLLARY 6.2:
(a) If $X$ is a noncritical \spsc\ extension of $A\subset X$ then every 
\holo\ submersion $f\colon A\to \C^q$ $(q<\dim X)$ can be 
approximated uniformly on $A$ by \holo\ submersions $\wt f\colon X\to \C^q$.
\item{(b)} Let $\Omega\subset \C^n$ be a convex open set.
Any \holo\ submersion $f\colon \Omega\to\C^q$ $(q<n)$ can be approximated 
uniformly on compacts by  submersions $\C^n\to\C^q$.

\demo Proof of Proposition 6.1:
Let $z=(z_1,\ldots,z_n)=(x_1+iy_1,\ldots, x_n+iy_n)$ denote the coordinates 
on $\C^n$. Let
$$
	P=\{z\in \C^n\colon |x_j|<1,\ |y_j|<1,\ j=1,\ldots,n\}
$$
denote the open unit cube. Set $P'=\{z\in P\colon y_n=0\}$.
 
Let $A, B\subset X$ be compact sets in $X$. We say that 
{\it $B$ is a convex bump on $A$} if there exist an open set
$U\subset X$ containing $B$, a biholomorphic map $\phi\colon U\to P$ 
onto $P\subset\C^n$, and smooth strongly concave functions 
$h, \wt h \colon P'\to [-a,a]$ for some $a<1$ such that 
$h\le \wt h$, $h=\wt h$ near the boundary of $P'$, and
$$
   \eqalign{ \phi(A\cap U) &= \{z\in P\colon 
         y_n \le h(z_1,\ldots,z_{n-1},x_n)\}, \cr
   \phi((A\cup B)\cap U) 
   	&= \{z\in P \colon  
  	 y_n \le \wt h  (z_1,\ldots,z_{n-1},x_n)\}. \cr}
$$

Suppose now that $A\subset \wt A$ is a noncritical \spsc\ extension in $X$.
By an elementary geometric argument,  using Narasimhan's lemma 
on local convexi\-fi\-ca\-tion of \spsc\ domains, there is a finite sequence 
$A=A_0\subset A_1\subset \ldots\subset A_{k_0}=\wt A$ of compact \spsc\ 
domains in $X$ such that for every $k=0,1,\ldots,k_0-1$ we have
$A_{k+1}=A_k\cup B_k$, where $B_k$ is a convex bump on $A_k$ as defined above. 
(For details see Lemma 12.3 in [HL2]. Similar `bumping constructions' 
had been introduced by Grauert and were used in the Oka-Grauert theory; 
see [Gro4], [HL3], [FP1], [FP2], [FP3].) Hence Proposition 6.1 follows 
immediately from the following. 

\proclaim LEMMA 6.3: Assume that $X$ is a Stein manifold, $A\subset X$
is a smooth compact \spsc\ domain, and $B\subset X$ is a convex bump on $A$.
Given a \holo\ submersion  $f \colon A\to \C^q$ $(q<\dim X)$, 
there exists for every $\e>0$ a \holo\ submersion $\wt f\colon A\cup B \to \C^q$ 
satisfying $\sup_{A\cap B} |f-\wt f|<\e$. If $X_0\subset X$ is a closed
complex subvariety such that $X_0\cap B=\emptyset$, we can choose
$\wt f$ such that it agrees with $f$ to a given finite order along 
$X_0 \cap A$.

\demo Proof:
We use the notation introduced above. 
Recall that $h$ and $\wt h$ have range in $[-a,a]$ for some $a<1$. 
Choose $c\in (a,1)$ sufficiently close to $1$ such that the 
(compact) support of $\wt h-h$ is contained in $cP'$. Let 
$L\colon=c\bar P \subset \C^n$ and $\wt L\colon=\phi^{-1}(L) \subset U$. 
Increasing $c<1$ towards $1$ we may assume that $B\subset \wt L$. 
Set $\wt K = A \cap \wt L$ and $K=\phi(\wt K) \subset P$. The pair of compact 
sets $K, L \subset \C^n$ satisfies the hypothesis of Proposition 3.3 
with respect to the splitting $z=(z',z'')\in\C^n$, with
$z'=(z_1,\ldots,z_{n-2}) \in\C^{n-2}$ and $z''=(z_{n-1},z_n)\in \C^2$.
Applying Proposition 3.3 (with $r=n-2$, $s=2$) we obtain a \holo\ 
submersion $g$ from a \nbd\ of $\wt L$ to $\C^q$ which approximates
$f$ uniformly in a \nbd\ of $\wt K$. Since $B \subset \wt L$ and 
$A\cap B\subset A\cap \wt L = \wt K$, $g$ is defined in a \nbd\ of $B$ 
and it approximates $f$ uniformly in a \nbd\ of $A\cap B$.  
By Lemma 5.1 we have $f=g\circ \g$ for a biholomorphic map $\g$ close to 
the identity in a \nbd\ of $A\cap B$ in $X$. Splitting $\g=\b\circ \a^{-1}$ 
by Theorem 4.1 we obtain $f\circ\a=g\circ \b$ in a \nbd\ of $A\cap B$, 
and hence the two sides define a \holo\ submersion 
$\wt f \colon A\cup B \to\C^q$. The same proof applies
with interpolation on $X_0$.
\medskip

%
%
%
%
In the remainder of this section we treat the critical case.
Let $p$ be a critical point of $\rho$, with Morse index $k$.
If $k=0$ then $\rho$ has a local minimum at $p$, and a new connected 
component appears in $\{\rho<c\}$ as $c$ passes $\rho(p)$. We can 
trivially extend $f$ to this new component by taking any 
local submersions to $\C^q$ near $p$. 
In the sequel we only treat the case $k\ge 1$. It is no loss of generality
to assume $\rho(p)=0$. Choose $c_0>0$ such that $p$ is the only critical 
point of $\rho$ in $[-c_0,3c_0]$.  In the following three subsections 
we explain how to approximately extend a submersion $f$ from
$\{\rho\le -c_0\}$ to $\{\rho\le +c_0\}$.

%
%
%
%
\medskip
{\it 2. Smooth extension across a handle.} 
Recall that $k\in \{1,\ldots, n\}$ is the index of $p$.
Write $z=(z',z'') = (x'+iy',x''+iy'')$, where $z'\in \C^k$ and $z''\in \C^{n-k}$.
Denote by $P\subset \C^n$ the open unit polydisc. By Lemma 3 in [HW2, p.\ 166] 
(see also Lemma 2.5 in [HL]) there is a \nbd\ $U\subset X$ of $p$ and a 
biholomorphic coordinate map $\phi\colon U\to P$, with $\phi(p)=0$, such that 
the function $\wt \rho(z) := \rho(\phi^{-1}(z))$ is given by
$$
	\wt \rho(z) = Q(y',z'')-|x'|^2,\quad 
	Q(y',z'')=\langle Ay',y'\rangle + \langle By'',y''\rangle + |x''|^2. \eqno(6.1)
$$
Here $\langle \cdotp,\cdotp\rangle$ is the Euclidean inner product and 
$A$, $B$ are positive definite symmetric matrices such that all eigenvalues 
of $A$ are larger than $1$ (thus $A>I$ and $B>0$). Furthermore one may
diagonalize $A$ and $B$. 

We may assume that $c_0<1$. Choose $c\in (0,c_0)$. By the noncritical case
we may assume that $f$ has already been extended to $\{\rho < -c/2 \}$. 
The set $E\subset U$ defined by 
$$
	\phi(E) = \{(x'+iy',z'')\colon  y'=0,\ z''=0, |x'|^2\le c \}
				    		                     \eqno(6.2) 
$$
is a $k$-dimensional handle attached from the outside to 
$\{\rho\le -c\}$ along the $(k-1)$-sphere $bE\subset \{\rho= -c\}$. 

In a \nbd\ of $E$ we may consider $f$ as a function of $z$. 
We identify $x\in \R^n$ with $x+i0\in\C^n$. The components 
$\theta_j$ of the $q$-coframe $\theta$ are expressed in the $z$ coordinates 
by  $\theta_j(z) = \sum_{l=1}^n \theta_{j,l}(z)\, dz_l$ where 
$\theta_{j,l}$ are continuous functions and the $q\times n$ matrix 
$\wt J=\pmatrix{\theta_{j,l}}$ has maximal complex rank $q$ at each point.
For $x\in E$ near $bE$ we have 
$\theta_{j,l}(x)=\di f_j/\di z_l(x)= \di f_j/\di x_l(x)$. 

Denote by $M_{q,n}\simeq \C^{q\times n}$ the set of 
all complex $q\times n$ matrices and let $M^*_{q,n}$ consist of all 
matrices of rank $q$ in $M_{q,n}$. 

\proclaim LEMMA 6.4: There is a $c'\in (0,c)$ such that $f$ and all its 
partial derivatives $\di f/\di z_l$  extend smoothly to $\{\rho\le -c'\}\cup E$
(without changing their values on $\{\rho\le -c'\}$) such that the Jacobian 
matrix $J(f)=(\di f_j/\di z_l)$ of the extension has complex rank $q$ 
at each point of $E$, and $J(f)$ can be connected to  
$\wt J=\pmatrix{\theta_{j,l}}$ by a homotopy of maps to $M^*_{q,n}$ 
which is fixed on $\{\rho\le -c'\}\cap E$.

Lemma 6.4 is obtained from a {\it convex integration lemma} 
due to Gromov [Gro2, Lemma 3.1.3]. We state the special case 
which is needed. Fix numbers $0<r<R$, $\d>0$, and let  
$$
	D = \{ x\in \R^n \colon |x'|\le R,\ |x''|\le \d\}, \quad
        A = \{ x\in \R^n \colon r\le |x'|\le R, |x''|\le \d\}. 
$$
%
%
%
%
\proclaim LEMMA 6.5: 
Assume that $f=(f_1,\ldots,f_q)\colon A\to \C^q$ $(q\le n)$ is a smooth map 
whose Jacobian $J(f) = \pmatrix{\di f_j/\di x_l}$ has complex rank $q$ 
at each point. If there exists a continuous 
map $\wt J\colon D\to M^*_{q,n}$ with $\wt J|_A=J(f)$ then there 
is a smooth map $\wt f \colon D\to \C^q$ such that (i) $\wt f|_A=f$,
(ii) the Jacobian $J(\wt f)$ has range in $M^*_{q,n}$, and (iii)
$J(\wt f)$ is homotopic to $\wt J$ through maps $D\to M^*_{q,n}$ 
which are fixed on $A$. If $q\le n-[{k\over 2}]$ then such $\wt J$ and 
$\wt f$ always exist.

{\it Proof.}
We have $M^*_{q,n}=M_{q,n}\bs \Sigma$ where the $\Sigma$
consist of all matrices of rank less than $q$. We claim that $\Sigma$ 
is an algebraic subvariety of complex codimension $n-q+1$ in 
$M_{q,n}\simeq \C^{q\times n}$. Assume that $B\in\Sigma$
has rank $q-1$. Choose $1\le j_1 < j_2<\ldots < j_{q-1}\le n$
such that the corresponding columns of $B$ are linearly independent.
Locally near $B$ the set $\Sigma$ is defined by vanishing of the 
determinants obtained by adding to the columns $j_1,\ldots,j_{q-1}$
any of the remaining $n-q+1$ columns of $B$. Locally this gives $n-q+1$ 
independent polynomial equations for $\Sigma$. A similar
argument holds when $B$ has rank $<q-1$. (See also Proposition 2 in [Fo2].)

We are looking for an extension $\wt f \colon D\to \C^q$ of $f$ whose 
Jacobian $J(\wt f)$ misses $\Sigma$. If $k< 2(n-q+1)$ (which is 
equivalent to $q\le n-[{k\over 2}]$) Thom's jet transversality theorem 
([Tho] or [GG, p.\ 54]) gives a maximal rank extension of $f$ and its 
full one-jet from $A$ to the $k$-dimensional disc 
$D_k=\{(x',0)\colon |x'|\le R\}$, and hence to an open \nbd\ 
$V \subset \R^n$ of $A\cup D_k$.  Clearly there exists 
a diffeomorphism $\psi\colon D\to \psi(D)\subset V$ which equals 
the identity on $A$. Then $\wt f= f\circ \psi$ has the desired properties.

The general case of Lemma 6.5 follows from 
Gromov's {\it convex integration lemma} [Gro2, Lemma 3.1.3].
(This can also be found in Section 2.4 of [Gro3]; see especially 
(D) and (E) in [Gro3, 2.4.1.]. Another source is Sect.\ 18.2 of [EM]; 
see especially Corollary 18.2.2.) 
To apply Gromov's lemma we consider $M_{q,n}$ as the space of all one-jets of 
smooth maps $D\to \C^q$ at any point $x\in D$ (that is, the space 
of all first order partial derivatives at $x$, ignoring the image point).
The open set $\Omega= M^*_{q,n} \subset M_{q,n}$ defines a differential 
relation of order one which is {\it ample in the coordinate directions} 
(see [Gro2] or Sect.\ 18.1\ in [EM] for a definition of this notion),
and the stated results follows from the convex integration lemma.

Ampleness of $\Omega$ in the coordinate directions means the following. 
Choose $l\in \{1,\ldots,n\}$ and fix in an arbitrary way the entries 
of a $q\times n$ matrix which do not belong to the column $l$
(these represent the partial derivatives $\di f_j/\di x_k$
for $k\ne l$ at some point $x$). Let $\Omega'\subset \C^q$ consist of all 
vectors whose insertion in the $l$-th column gives a matrix of maximal 
rank $q$ (thus belonging to $\Omega$). $\Omega$ is ample in the
coordinate directions if every such set $\Omega'$ is either empty 
or else the convex hull of each of its connected components equals $\C^q$. 
In our case $\Omega'$ is either empty, the complement of a 
complex hyperplane in $\C^q$, or all of $\C^q$, depending on the rank of the 
initial $q\times (n-1)$ matrix. This completes the proof of Lemma 6.5.

\demo Proof of Lemma 6.4: Let $A\subset D$ be subsets of $U\subset X$ defined by 
$$\eqalign{
	\phi(D) &= \{(x'+i0',x''+i0'') \colon |x'|^2\le c,\ |x''|\le \d\}, \cr
        \phi(A) &= \{ z\in \phi(D) \colon r \le |x'|^2 \le c\}.
        } 
$$
Choosing $\d>0$ sufficiently small and $r < c$ sufficiently close to $c$ 
we insure that $A \subset \{\rho< -c/2\}$, and hence $f|_A \colon A\to\C^q$ is 
a well defined smooth map with differential of maximal complex rank $q$. 
Lemma 6.5 gives the desired smooth extension to $D$ as well as a homotopy 
of $q$-coframes which is fixed on $A$.  If $c'<c$ is chosen sufficiently close 
to $c$ then $D\cap \{\rho \le -c'\} \subset A$ and hence Lemma 6.4 holds 
for such $c'$.

%
%
\medskip
{\it 3. Holomorphic approximation.} 
Let $f$ be given by Lemma 6.4.  In this subsection we prove

\proclaim LEMMA 6.6: For every $\eta >0$ there exist an open 
\nbd\ $\Omega\subset X$ of the set $K=\{\rho \le -c\}\cup E$  
and a holomorphic submersion $\wt f\colon \Omega\to\C^q$ 
such that $|f-\wt f|_K <\eta$, $|df-d\wt f|_E< \eta$,
and $d\wt f$ is $q$-coframe homotopic to $\theta$.

Here $|f|_K$ is the uniform norm of $f$ on $K$, and $|df|_E$ is the norm
of its differential on $E$, measured in a fixed Hermitean 
metric on $TX$.

\demo Proof:  We need an improved version of Theorem 4.1 from [HWe]. 
We first show that $K$ is $\cO(X)$-convex and hence admits a basis of Stein \nbd s.
We use the notation from Subsect.\ 2. Choose $L\subset U$ 
such that $\phi(L)= r\bar P$ for some $r<1$ very close to $1$. Then 
$\phi(K\cap L) = \{z\in r\bar P\colon \wt\rho(z) \le -c\} \cup \phi(E)$.
Clearly each of the sets $\{\wt \rho\le -c\} \cap r\bar P$ and $\phi(E)$ is 
\pc\ in $\C^n$. The holomorphic polynomial $h(z)=z_1^2+\ldots +z_k^2$ maps 
$\phi(E)$ to the interval $[0,c]$, $\phi(bE)$ to the point $c$, and from (6.1) 
we easily see that $\Re h > c$ on $\{\wt \rho\le -c\}\bs \phi(E)$. 
Thus $h$ separates the two sets and hence their union is polynomially 
convex (Lemma 29.21 in [Sto]). $\cO(X)$-convexity of $\{\rho\le -c\}\cup E$ 
follows by a usual patching argument, using strong plurisubharmonicity 
of $\rho$ (see Lemma 1 in [Ro]). 

Choose a constant $\wt c \in (c',c)$. By Lemma 4.3 in [HWe] there is a smooth 
map $g \colon X\to \C^q$ satisfying

\smallskip
\item{(i)}  $g=f$ on $\{\rho\le - \wt c \} \cup E$,  
\item{(ii)} $dg_x=df_x$ for each $x\in E$,
\item{(iii)} $g$ is $\dibar$-flat on $E$, i.e., 
$D^r(\dibar g)|_E=0$ for all $r\in\N$. 

\smallskip
Here $D^r$ denotes the total derivative of order $r$. The cited lemma 
is proved in [HWe] for $X=\C^n$, but the result is local and holds 
for any smooth totally submanifold $E$ in a complex manifold. 
(One may use partitions of unity along $E$ which are $\dibar$-flat on $E$; 
see Lemma 2.3 in [FL\O].) If $E$ is of class $\cC^m$ then (iii) holds for 
$r\le m-1$.

Fix an integer $m\ge n+1$. Let $\Omega_\e=\{x\in X\colon d(x,K)<\e\}$.
In the proof of Theorem 4.1 in [HWe] on pp.\ 15-16 the authors obtained 
for each sufficiently small $\e>0$ a map $w_\e\colon \Omega_e\to\C^q$ 
satisfying $\dibar w_\e=\dibar g$ in $\Omega_\e$ 
and $||w_\e||_{L^2(\Omega_\e)} = o(\e^m)$ as $\e\to 0$. 
(The proof in [HWe] remains valid in any Stein manifold by applying 
the appropriate $\dibar$-results from [H\"o1].) 
On $\Omega_{\e/2}$ this gives a uniform estimate $|w_\e|=o(\e^{m-n})$ 
[HWe, p.\ 16] as well as $|D^r w_\e|=o(\e^{m-n-r})$ 
(Lemma 3.2 in [FL\O]). By construction the map 
$f_\e=g-w_\e \colon \Omega_\e \to\C^q$ is holomorphic and satisfies 
$|f_\e-f|=o(\e)$, $|df_\e-df|=o(1)$ on $\Omega_{\e/2}$ as $\e\to 0$. 
Hence for sufficiently small $\e>0$ the map $f_\e$ is a \holo\ submersion in 
an open \nbd\ $\Omega$ of $K$, with $df_\e$ close to $df$ and hence 
homotopic to $\theta$. This proves Lemma 6.6.

\medskip
{\it Remark.}  More precise approximation results on totally real 
submanifolds have been obtained by integral kernels; 
see [HW1], [RS] and [FL\O]. The paper [FL\O] contains optimal results 
on approximation of $\dibar$-flat functions in tubes around totally 
real submanifolds.

%
%
\medskip
{\it 4. Extension across the critical level.}
Let $\wt f \colon \Omega \to\C^q$ be the submersion furnished by
Lemma 6.6. To simplify the notation we drop the tilde.
The purpose of this subsection is to approximately extend $f$ 
across the critical level $\{\rho=0\}$ by applying the noncritical case 
(Proposition 6.1) with a different \spsh\ function $\tau$ given 
by Lemma 6.7 below. Once this is done, we switch back to $\rho$ 
(perhaps sacrificing some of the gained territory) and continue
(by the noncritical case) to its next critical level. 

We shall use the notation established in Subsect.\ 6.1.
Let $\phi\colon U\to  P \subset \C^n$ be a coordinate map 
as in the proof of Proposition 6.1 such that 
$\wt \rho= \rho\circ\phi^{-1}$ is given by  (6.1).
Let $c_0>0$ be the constant chosen in the paragraph preceding 
Subsect.\ 6.2. By the noncritical case we may decrease $c_0$ to insure   
$$
	\{(x'+iy',z'') \in\C^n \colon |x'|^2 \le c_0,\ Q(y',z'')\le 4c_0\} \subset  P.
$$
Denote by $E$ the handle (6.2) but with $c$ replaced by $c_0$ (thus $\rho=-c_0$ on $bE$). 
Let $\l_1> 1$ denote the smallest eigenvalue of the matrix $A$. Choose a number 
$1 < \mu < \l_1$ and set $t_0 = (1-1/\mu)^2 c_0$.

\proclaim LEMMA 6.7: 
There exists a smooth \spsh\ function $\tau$ on $\{\rho< 3c_0\} \subset X$
which has no critical values in $(0,3c_0) \subset \R$ and satisfies 
\item{(i)} $\{\rho\le -c_0\} \cup E \subset 
\{\tau\le 0\} \subset \{\rho\le -t_0\}\cup E$, and
\item{(ii)} $\{\rho \le c_0\} \subset \{\tau \le 2c_0\} \subset \{\rho< 3c_0\}$.
\smallskip

Using Lemma 6.7 we complete the crossing of the critical level $\{\rho=0\}$ 
as follows. Consider the family of sublevel sets $\{\tau \le c\}$ as $c$ 
increases from $0$ to $2c_0$. Each of them contains the handlebody 
$\{\rho\le -c_0\} \cup E$. Property (i) implies that for sufficiently small 
$c>0$ we have $\{\tau\le c\} \subset \Omega$. 
By Proposition 6.1 (the noncritical case) we can approximate 
$f$ uniformly on $\{\tau \le c\}$ by a submersion $\wt f$ defined in a
\nbd\ of $\{\tau \le 2c_0\}$. By (ii) $\wt f$ is defined on 
$\{\rho\le c_0\}$ and $d\wt f$ is $q$-coframe homotopic to $\theta$.
Since $c_0>0$, this completes the extension across the critical level $\{\rho=0\}$.
Hence Theorems 2.5 and 2.6 are proved except for the interpolation 
on a subvariety (Subsect.\ 5).
\medskip

In the proof of Lemma 6.7 we shall need a criterion for strong 
plurisubharmonicity of certain functions modeled on (6.1).

\proclaim LEMMA 6.8: 
Let $A>0$ be a symmetric real $n\times n$ matrix with the smallest 
eigenvalue $\l_1>0$. If a $\cC^2$ function 
$h\colon I\subset \R_+ \to\R$ satisfies 
$$ 
	\dot h <\l_1 \quad{\rm and}\quad 
	2t \ddot h + \dot h < \l_1       \qquad(t\in I)         \eqno(6.3)
$$
then the function $\tau(z)=\langle Ay,y\rangle - h(|x|^2)$
is \spsh\ on $\{z=x+iy \in\C^n\colon |x|^2\in I\}$.

\demo Proof:
Let $A=\pmatrix{a_{jl}}$. A calculation gives
$$ 
	\eqalign{ 
	-\tau_{z_j} &= x_j\dot h + i\sum_{s=1}^n a_{js}y_s \cr
	-2\tau_{z_j \bar z_l} &= 
	\cases{ 2x_j^2 \ddot h + \dot h -a_{jj},   & if $j=l$;\cr
	        2x_j x_l \ddot h  - a_{jl},        & if $j\ne l$. \cr} 
	    }
$$
Thus the complex Hessian 
$H_\tau=\pmatrix{ {\di^2 \tau \over \di z_j \di{\bar z}_l}}$ of $\tau$ satisfies
$$ 
	-2H_\tau = 2\ddot h \cdotp xx^t + \dot h I - A       
$$
where $xx^t$ is the matrix product of the column $x\in\R^k$ with the row $x^t$ 
and $I$ denotes the identity matrix. For any $v\in\R^n$ we have 
$
	\langle (x x^t) v,v \rangle = v^t x x^t v 
	= |\langle x,v\rangle|^2 
$
which lies between $0$ and $|x|^2|v|^2$. Hence $0\le xx^t \le |x|^2 I$. 
(Here we write $A\le B$ if $B-A$ is nonnegative definite.) 
At points $|x|^2=t$ where $\ddot h(t)\ge 0$ we thus get  
$-2H_\tau \le (2t\ddot h + \dot h)I - A < \l_1 I-A \le 0$ 
and hence $H_\tau>0$ (we used the second inequality in (6.3)). 
At points where $\ddot h<0$ we can omit $2\ddot h xx^t \le 0$
to get $-2H_\tau \le \dot h I -A \le (\dot h -\l_1)I <0$, so $H_\tau>0$.
Thus $H_\tau$ is positive definite which proves Lemma 6.8.

\demo Proof of Lemma 6.7: Recall that $1<\mu <\l_1$ and 
$t_0 = (1-1/\mu)^2 c_0$. We shall find a smooth convex increasing
function $h\colon \R \to [0,+\infty)$ satisfying  
\item{(i)}   $h(t)=0$ for $t\le t_0$, 
\item{(ii)}  $h(t)=t - t_1$ for $t\ge c_0$, where $t_1=c_0 - h(c_0) \in (t_0,c_0)$, and
\item{(iii)} for all $t\ge t_0$ we have $0\le \dot h \le 1$, 
$2t\ddot h + \dot h<\l_1$, and $t-t_1 \le h(t) \le t - t_0$.

We first consider the function   
$$
	\xi(t)= \cases{ 0,         & if $t \le t_0$; \cr
                 \mu \bigl( \sqrt t -\sqrt {t_0}\, \bigr)^2, 
                                   & if $t_0 \le t\le c_0$, \cr
                 t-c_0(1-1/\mu),   & if $c_0 \le t$.  \cr}
$$
On $[t_0,c_0]$ $\xi$ solves the initial value problem
$2t\,\ddot{\xi} + \dot \xi =\mu$, $\xi(t_0)=\dot \xi(t_0)=0$. 
It is $\cC^1$ and piecewise $\cC^2$, with discontinuities 
of $\ddot \xi$ at $t_0$ and $c_0$. The value of $t_0$ is chosen such that 
$\dot \xi(c_0)=1$. We have $\ddot \xi(t) = \mu\sqrt{t_0}/2\sqrt{t}^3 >0$
for $t\in [t_0,c_0]$, $\ddot \xi(t) = 0$ for $t$ outside this interval,
and $\int_{t_0}^{c_0} \ddot \xi(t) dt = 1$.

Choose a smooth function $\chi\ge 0$ which vanishes outside $[t_0, c_0]$, 
equals $\ddot \xi +\e$ on $[t_0+\d,c_0-\d]$ for small $\e,\d>0$, and  
interpolates  between $0$ and $\ddot \xi$ on the intervals
$[t_0,t_0+\d]$ and $[c_0-\d,c_0]$. We can find $\d,\e>0$
arbitrarily small such that 
$\int_{t_0}^{c_0} \chi(t)dt = \int_{t_0}^{c_0} \ddot \xi(t) dt= 1$. 
The function $h\colon \R_+\to\R_+$ obtained by integrating $\chi$ 
twice with the initial conditions $h(t_0)=\dot h(t_0)=0$ will satisfy 
the properties (i)--(iii) provided that $\e$ and $\d$ were chosen
sufficiently small (since $h$ is then $\cC^1$-close to $\xi$ 
and $\ddot h \le \ddot \xi +\e$). In particular, 
$t_1=c_0-h(c_0) \approx c_0-\xi(c_0)=(1-1/\mu)c_0$
and hence $t_0< t_1 < c_0$. 

By Lemma 6.8 the function 
$$
	\wt \tau(z)= \langle Ay',y'\rangle - h(|x'|^2) 
	                   + \langle By'',y''\rangle + |x''|^2
	           =Q(y',z'') - h(|x'|^2) 
$$
is \spsh\ on $\C^n$. Recall that $\wt\rho(z)=Q(y',z'') - |x'|^2$.
The properties of $h$ imply 

\smallskip
\item{(a)} $\wt \rho\le \wt \tau \le \wt \rho+t_1$, 
\item{(b)} $\wt \rho+ t_0\le \wt \tau$ on the set $\{|x'|^2\ge t_0\}$, and
\item{(c)} $\wt \tau = \wt\rho + t_1$ on $\{|x'|^2\ge c_0\}$.

\smallskip
Let $V= \{\rho<3c_0\}\subset X$. We define $\tau \colon V \to\R$ by 
$\tau=\wt\tau\circ \phi$ on $U\cap V$ and $\tau=\rho + t_1$ on $V\bs U$. 
Property (c) implies that both definitions agree on $U\cap V\cap \{|x'|^2\ge c_0\}$ 
and hence $\tau$ is \spsh. The stated properties of $\tau$ follow immediately.
This completes the proof of Lemma 6.7.

%
%
\medskip
{\it 5. Interpolation along a complex subvariety.}
In this subsection we prove the interpolation statement in Theorem 2.5.
Recall the situation: 
\item{--} $X_0$ is a closed complex subvariety of a Stein manifold $X$, 
\item{--} $K$ is a compact $\cO(X)$-convex subset of $X$,
\item{--} $U\subset X$ is an open set containing $K\cup X_0$,
\item{--} $f\colon U\to\C^q$ is a holomorphic submersion such that
$df=\theta|_U$ for some $q$-coframe $\theta$ defined on $X$. 
\smallskip

Let $c$ be a regular value of $\rho$ such that $L=\{\rho\le c\}$ 
contains $K$ in its interior.  Our task is to find a \holo\ 
submersion $\wt f$ from an open \nbd\ of $L\cup X_0$ to $\C^q$ which 
approximates $f$ uniformly on $K$, it interpolates $f$ along $X_0$ 
to order $r\in \N$, and $d\wt f$ is $q$-coframe homotopic to $\theta$. 
The desired submersion $X\to\C^q$ is then obtained by a usual 
limiting process. For convenience of notation we take $c=0$ 
and $L=\{\rho\le 0\}$.

The set $K' := (K\cup X_0)\cap  \{\rho\le 1\}$ is $\cO(X)$-convex
and hence there exists a smooth \spsh\ exhaustion function 
$\tau \colon X\to \R$ such that $\tau<0$ on $K'$ and $\tau >0$ on $X\bs U$. 
We may assume that $0$ is a regular value of $\tau$ and the hypersurfaces 
$\{\rho=0\}$ and $\{\tau=0\}$ intersect transversely. The set  
$D_0=\{\tau \le 0\}$ is a smooth \spsc\ domain contained in 
the domain $U$ of $f$. The following lemma provides the main step.

\proclaim LEMMA 6.9:
For each $\e>0$ there exists a \holo\ submersion $g\colon \wt L\to \C^q$  
in an open set $\wt L\supset L$ such that $|g-f|<\e$ on $D_0\cap L$ 
and $g - f$ vanishes to order $r$ on  $X_0\cap \wt L$.

Assuming Lemma 6.9 we complete the proof of Theorem 2.5 
as follows. Cartan's theory gives $f'\in\cO(X)^q$ such that $f'-f$ 
vanishes to order $r$ on $X_0$, and finitely many
functions $\xi_j\in\cO(X)$ $(j=1,2,\ldots,m)$ which vanish to order $r$ 
on $X_0$ and generate the corresponding sheaf of ideals 
$\cJ^r_{X_0}$ on $L$ (but not necessarily on $X$). 
Since $g-f'\in\cO(L)^q$ vanishes to order $r$ on $X_0\cap L$, we have 
$g=f' + \sum_{j=1}^m \xi_j g_j$ for some $g_j\in \cO(L)^q$. Since $L$ is 
$\cO(X)$-convex, we can approximate each $g_j$ uniformly on a \nbd\ of $L$ 
by $\wt g_j\in\cO(X)^q$. The map  $\wt f=f'+ \sum_{j=1}^m \xi_j \wt g_j \colon X\to\C^q$ 
is holomorphic, $|\wt f- g|$ is small on a \nbd\ of $L$ (hence $|\wt f-f|$ is small 
on $D_0\cap L$), and $\wt f-f$ vanishes to order $r$ along $X_0$. 
If the approximations are sufficiently close then $\wt f$ is a submersion 
in a \nbd\ of $L\cup X_0$. This completes the induction step.

\demo Proof of Lemma 6.9:
Set $\rho_t=\tau + t(\rho-\tau)=(1-t)\tau + t\rho$ and let
$$
	D_t=\{\rho_t \le 0\} =\{\tau \le t(\tau-\rho) \} \qquad (t \in [0,1]).
$$
We have $D_0=\{\tau\le 0\}$, $D_1=\{\rho\le 0\}=L$, and $D_0\cap D_1\subset D_t$
for all $t\in [0,1]$. Let $\Omega= \{\rho<0,\ \tau>0\} \subset D_1\bs D_0$, 
$\Omega'=\{\rho>0,\ \tau<0\} \subset D_0\bs D_1$. Since $\tau-\rho>0$ in $\Omega$ 
and $\tau-\rho<0$ in $\Omega'$ it follows that, as $t$ increases from $0$ to $1$, 
the sets $D_t \cap L$ monotonically increase to $D_1=L$ while $D_t\bs L \subset D_0$ 
decrease to $\emptyset$. All hypersurfaces $\{\rho_t=0\}=bD_t$ intersect along the 
real codimension two submanifold $S= \{\rho=0\}\cap \{\tau=0\}$. Since 
$d\rho_t=(1-t)d\tau+t d\rho$ and the differentials $d\tau$, $d\rho$ are linearly 
independent along $S$, each hypersurface $bD_t$ is smooth near $S$. Since $\rho_t$ 
is a convex linear combination of \spsh\ functions, it is itself \spsh\ and hence 
$D_t$ is \spsc\ at every smooth point of $bD_t$. 

We investigate more closely the non-smooth points of $bD_t =\{\rho_t=0\}$ inside 
$\Omega$. The defining equation of $D_t\cap \Omega$ can be written as 
$\tau \le t(\tau-\rho)$ and, after dividing by $\tau-\rho>0$, as
$$   
	D_t\cap \Omega = \{x\in \Omega \colon h(x)= 
	{\tau(x)\over \tau(x)-\rho(x)} \le  t\}. 
$$
The equation $dh=0$ for critical points is equivalent to 
$(\tau -\rho)d\tau - \tau(d\tau -d\rho)= \tau d\rho-\rho d\tau =0$.
A generic choice of $\rho$ and $\tau$ insures that there are at most
finitely many solutions $p_1,\ldots,p_m\in \Omega$ and no solution
on $b\Omega$. A calculation shows that at each critical point 
the complex Hessians satisfy  $(\tau -\rho)^2 H_h = \tau H_\rho - \rho H_\tau$.
Since $\tau>0$ and $-\rho>0$ on $\Omega$, we conclude that $H_h>0$ at such points.
By a small modification of $h$ near each $p_j$ we can therefore assume 
that it is of the form (6.1) in some local holomorphic coordinates. 

If $c\in [0,1)$ is a regular value of $h|_\Omega$ then for $c'>c$ sufficiently 
close to $c$ (depending only on $h$) the domain $D_{c'}$ can be 
obtained from $D_c\cap D_{c'}$ by finitely many attachings of 
{\it convex bumps} (Subsect.\ 1). Indeed, the boundaries $bD_c$ and 
$bD_{c'}$ intersect transversely at very small angles along $S$
and are locally convexifiable. We begin by attaching small convex 
bumps to $D_c\cap D_{c'}$ along $S$ in order to enlarge $D_c \cap L$ 
to $D_{c'} \cap L$ locally near $S$ while keeping unchanged the part of the
set outside of $L$ (which equals $D_{c'}\bs L$). Each of the bumps 
may be chosen disjoint from $X_0$ and with finitely many bumps 
we can reach $D_{c'}$. By Lemma 6.3 every submersion 
defined in a \nbd\ of $D_c$ can be approximated uniformly on $D_c\cap D_{c'}$ 
by a submersion defined in a \nbd\ of $D_{c'}$ such that the two maps agree 
to order $r$ on $X_0$. (We use the interpolation version of Theorem 4.1.) 
If $0\le c_0 < c_1\le 1$ 
are such that $h|_\Omega$ has no critical values in $[c_0,c_1]$, we 
can subdivide $[c_0,c_1]$ into finitely many subintervals on which 
the above procedure applies. This explains the noncritical case. 

We have seen that the (finitely many) critical points of $h|_\Omega$ 
are of the form (6.1) and hence the method developed in Subsections 2--4 
can be applied to cross every critical level of $h$. Lemma 6.6 with 
interpolation on $X_0$ (which does not intersect the handle $E$)
is a trivial addition. 

Together these two methods show that we can approximate a \holo\ submersion $f$, 
defined in a \nbd\ of $D_0$, uniformly on $L\cap D_0 \supset K$ by a 
submersion $g$ defined in a \nbd\ of $L$ such that $g-f$ vanishes 
to order $r$ on $X_0$. This completes the proof of Lemma 6.9.

\medskip
{\it Proof of Corollary 2.10.} 
The hypothesis implies that the normal bundle of $V$ in $X$ is 
isomorphic to $N|_V$ and hence is trivial. By the Docquier-Grauert 
theorem [DG] there exist functions $g_1,\ldots,g_q\in\cO(U)$ 
whose common zero set equals $V$ and whose differentials $dg_1,\ldots,dg_q$ are
linearly independent along $V$. If $dg$ extends to a $q$-coframe 
on $X$ then Theorem 2.5 furnishes a submersion $f\colon X\to\C^q$ 
such that $f-g$ vanishes to second order along $V$. This implies that 
$V$ is a union of connected components of $f^{-1}(0)$. 

In general we must replace $g_1,\ldots,g_q$ by a different set of 
defining functions for $V$ to insure the $q$-coframe extendability. 
Choose a complex subbundle $E \subset TX$ such that $TX=E \oplus N$ 
and $E=\ker dg$ in a \nbd\ of $V$ (in particular, $E|_V = TV$).  
Let $\Theta \subset T^*X$ be the conormal bundle with fibers 
$\Theta_x=\{\omega\in T_x^*X\colon \omega(v)=0\ {\rm for\ all\ } v\in E_x \}$. 
From $\Theta \simeq (TX/E)^* \simeq N^*$ we see
that $\Theta$ is trivial. Hence there exists a $q$-coframe 
$\theta=(\theta_1,\ldots,\theta_q)$ on $X$ which spans $\Theta$ 
and is holomorphic near $V$. By construction the differential $dg_1,\ldots,dg_q$ 
also spans $\Theta$ near $V$ and hence $\theta_j=\sum_{k=1}^q a_{jk} dg_k$
for some holomorphic functions $a_{jk}$ in a \nbd\ of $V$.
Set $h_j =\sum_{k=1}^q a_{jk} g_k$ for $j=1,\ldots,q$.
Then $dh_j=\theta_j$ at points of $V$ (since the term obtained by
differentiating $a_{jk}$ is multiplied by $g_k$ which vanishes on $V$).
Let $\chi$ be a smooth function on $X$ which equals one in a 
small \nbd\ of $V$ and equals zero outside of a slightly larger 
\nbd. If these neighborhood are chosen sufficiently small
then $\wt \theta = \chi dh + (1-\chi) \theta$ is a
$q$-coframe on $X$ which equals $dh$ near $V$. Hence we can 
apply Theorem 2.5 to $h$ as explained above.  

Assume now  $\dim V\le [{n\over 2}]$ so the rank of its (trivial) 
normal bundle $N_V$ is at most $[{n+1\over 2}]$. It suffices 
to show that $N_V$ extends to a trivial subbundle $N\subset TX$. To see this, 
recall that the pair $(X,V)$ is homotopy equivalent to a relative 
CW-complex of dimension at most $n$ [AF]. The standard topological method 
of extending  sections over cells gives the following: 
{\it If $E\to X$ is a complex vector bundle of rank $k>n/2$ then 
a nonvanishing section of $E$ over $V$ extends to a nonvanishing 
section of $E$ over $X$}. Indeed the obstruction to extending a section 
from the boundary of an $m$-cell to its interior lies in the homotopy 
group $\pi_{m-1}(S^{2k-1})$ which vanishes if $m<2k$. Our complex
only contains cells of dimension $\le n$ which gives the stated result.
Using this inductively we see that the linearly independent sections 
generating $N_V \subset TX|_V$ extend to linearly independent 
sections over $X$ generating a trivial subbundle $N\subset TX$.

\medskip
{\it Proof of Corollary 2.11.} 
By [DG] there exists an open set $U\subset X$ containing $V$ and a
\holo\ submersion $\pi\colon U\to V$ which retracts $U$ onto $V$. 
Choose a \holo\ subbundle
$H\subset TU$ such that $TU=H\oplus \ker d\pi$ and $H|_V=TV$.
The map $f^0 :=f\circ \pi \colon U\to \C^q$ is a \holo\ submersion
with $f^0|_V=f$. By the assumption there is a $q$-coframe 
$\theta=(\theta_1,\ldots,\theta_q)$ on $X$  
satisfying $\iota^* \theta_j=df_j$ for $j=1,\ldots, q$.
Choose a smooth cut-off function $\chi\colon X\to [0,1]$ with 
support in $U$ such that $\chi=1$ in a smaller open \nbd\ $U_1 \subset U$ 
of $V$. The $(1,0)$-forms $\wt \theta_j := \chi df^0_j +(1-\chi)\theta_j$
($j=1,\ldots,q$) are well defined on $X$ and are $\C$-linearly
independent, except perhaps on the set where $0<\chi<1$. However,
if we choose $\chi$ to be supported in a sufficiently thin \nbd\ of $V$ 
then these forms are also independent there since the $H$-components
of the $q$-coframes $\theta$ and $df^0$ agree on $H|_V=TV$, and hence
are close to each other over an open \nbd\ of $V$. It remains to 
apply Theorem 2.5 to obtain a submersion $F\colon X\to \C^q$ extending
$f$. If $q\le [{n+1\over 2}]$ then the $q$-coframe $df^0$ extends
from a small \nbd\ of $V$ to all of $X$ by the same argument as in 
the proof of Corollary 2.10, using the fact that the pair $(X,V)$
is homotopic to a relative CW-complex of dimension $\le \dim X$.

%
%
%
%
\beginsection 7. Holomorphic sections transverse to a foliation

A complex vector bundle $\pi\colon N\to X$ of rank $q$ admits 
{\it locally constant transition functions} if there is an open covering 
$\{U_i\}_{i\in\N}$ of $X$ and fiber preserving homeomorphisms 
$\phi_i\colon N|_{U_i}=\pi^{-1}(U_i) \to U_i\times \C^q$ 
with transition maps
$$
	\phi_{ij}(x,z)=\phi_i\circ\phi_j^{-1}(x,z) =
	\bigl( x,h_{ij}(z) \bigr)      
	\quad (x\in U_i\cap U_j,\ z\in \C^q)                     
$$ 
where $h_{ij}\in GL_q(\C)$ is independent of the base point $x\in U_i\cap U_j$.
The {\it structure group} $\Gamma\subset GL_q(\C)$ of $N$, generated by all $h_{ij}$'s,
is totally disconnected but not necessarily discrete. (Such $N$, also called
a {\it flat bundle}, is determined by a representation 
$\alpha \colon \pi_1(X)\to GL_q(\C)$; its pull-back to 
the universal covering $\wt X$ of $X$ is a trivial bundle over $\wt X$. 
This will not be used in the sequel.)

\proclaim THEOREM 7.1: 
Let $X$ be a Stein manifold. If $E$ is a complex subbundle of the tangent bundle
$TX$ such that $N=TX/E$ admits locally constant transition functions then 
$E$ is homotopic (through complex subbundles of $TX$)
to the tangent bundle of a nonsingular \holo\ foliation of $X$.

Theorem 7.1 extends Corollary 2.9 in which $N=TX/E$ was assumed 
to be trivial. The analogous result concerning smooth foliations on smooth open manifolds 
was proved by Gromov [Gro1] and Phillips ([Ph2], [Ph3], [Ph4]), and on closed manifolds 
by Thurston ([Th1], [Th2]). (See also [God, pp.\ 65-66] and [Gro3, p.\ 102].)
The smooth analogue of Theorem 7.1 applies to any smooth codimension one 
subbundle $E\subset TX$ (since any real line bundle admits a totally disconnected
structure group). On the other hand, a complex line bundle $N\to X$ over a Stein manifold 
admits such a structure group only if its first Chern class $c_1(N)\in H^2(X;Z)$
is a torsion element of this group.

\demo Proof of Theorem 7.1:
Since the transition functions $h_{ij}$ do not depend on the base point,
the product foliations over the sets $U_j \in \cU$ define a global 
\holo\ foliation $\cH$ of $N$ such that the zero section of $N$ is a 
union of leaves (one for each connected component of $X$). 
More precisely, if $U_i\cap U_j\ne\emptyset$ and $z\in\C^q$
then $\phi^{-1}_i(U_i \times\{h_{ij}(z)\})$ and $\phi^{-1}_j(U_j \times\{z\})$ 
belong to the same leaf of $\cH$. The tangent bundle of $N$ decomposes as 
$TN=H\oplus V$ where the {\it horizontal component\/} $H:=T\cH$ is the tangent 
bundle of $\cH$ and the {\it vertical component} $V$ is the tangent bundle of 
the foliation $N_x =\pi^{-1}(x)$ $(x\in X)$. Denote by $\tau \colon TN\to V$ 
the projection onto $V$ with kernel $H$. Observe that $V$ is just the pull-back 
of the vector bundle $N\to X$ to the total space by the projection map $\pi$,
and for every section $f\colon X\to N$ of $\pi$ we have $f^*V=N$.

If $f\colon X \to N$ is a holomorphic section transverse to $\cH$ 
(this requires $q\le n=\dim X$) then the intersections of $f(X) \subset N$ 
with the leaves of $\cH$ defines a \holo\ foliation $\cH_f$ of $X$, of 
dimension $k=n-q$, whose tangent bundle $T\cH_f \subset TX$ has fibers 
$(T\cH_f)_x = \{\xi \in T_x X \colon\ \tau \circ df_x(\xi) =0 \}$. 
Transversality of $f$ to $\cH$ means that the vector bundle map 
$$
	f' := f^*\circ\tau\circ df \colon TX\to f^*V=N		
$$ 
is surjective and hence induces an isomorphism of 
$TX/T\cH_f$ onto $N$. In particular, $N$ is the 
normal bundle of any such foliation $\cH_f$. 

To prove Theorem 7.1 we construct a holomorphic  section
$f\colon X\to N$ transverse to $\cH$ and a complex vector bundle 
injection $\iota\colon N\to TX$ (not necessarily holomorphic)
such that the subbundle $T\cH_f \subset TX$ is homotopic to $E$ 
and $f'\circ\iota\colon N\to N$ is a complex vector bundle automorphism 
homotopic to the identity through complex vector bundle automorphisms of $N$.

On every sufficiently small open set $U\subset X$ we have $N|_U\simeq U\times\C^q$ 
and the restriction of $\cH$ to $N_U$ has leaves $U\times\{z\}$ $(z\in\C^q)$.
Any such $U$ will be called {\it admissible}.
A section of $N$ over such $U$ is of the form $f(x)=(x,\wt f(x))$ where 
$\wt f\colon  U\to \C^q$, and $f$ is transverse to $\cH$ \iff\ $\wt f$ is 
a submersion to the fiber $\C^q$. This reduces every local problem 
in the construction of a transverse section to the corresponding problem for 
submersions. 

Choose a \spsh\ Morse exhaustion function $\rho\colon X\to\R$
and an initial embedding $\tau\colon N\to TX$ such that 
$TX=E\oplus \iota(N)$. Suppose $f$ is a transverse \holo\ section, 
defined on a sublevel set of $\rho$, such that $\ker (\tau\circ df)$ is complementary 
to $\iota(N)$  and $f'\circ\iota$ is homotopic to the identity over the domain 
of $f$. We inductively enlarge the domain of $f$ as in the proof of 
Theorem 2.5. Whenever we change $f$ the injection $\iota$ is changed 
accordingly (by a homotopy of injections $N\to TX$) such that $f'\circ\iota$ 
remains homotopic to the identity on $N$. We must explain the following two steps.

\smallskip
\item{(a)} Suppose that $(A,B)$ is a special Cartan pair in $X$ such that
$B$ is a convex bump on $A$ contained in an admissible set $U\subset X$
(Subsect.\ 6.1). Given a transverse section $f\colon \wt A\to N$ 
in a \nbd\ of $A$, find a transverse section $F$ in a \nbd\ of $A\cup B$ 
which approximates $f$ uniformly on $A$. (The homotopy conditions trivially extend 
from $A$ to $A\cup B$.) A solution to this problem will complete the proof
in the noncritical case (compare with Proposition 6.1).

\smallskip
\item{(b)} Extend a transverse section across a critical level of 
a $\rho$. At this step we shall need the homotopy condition on $f'\circ\iota$.

\smallskip
Part (a) is proved as in Proposition 6.1 with one minor change.
On $\wt A\cap U$ we have $f(x)=(x,\wt f(x))$ where $\wt f$ is 
a submersion to $\C^q$. We approximate $\wt f$ uniformly 
in a \nbd\ of $A\cap B$ by a submersion $\wt g\colon \wt B\to \C^q$ 
defined in a \nbd\ of $B$, find a transition map $\gamma$ such 
that $\wt f=\wt g\circ \gamma$ in a \nbd\ of $A\cap B$, 
and split $\gamma=\beta\circ\alpha^{-1}$ by Theorem 4.1. This gives 
$\wt f\circ \alpha = \wt g\circ \beta$ in a \nbd\ of $A\cap B$
which defines a transverse \holo\ section $F$ in a \nbd\ of 
$(A\cup B)\cap U$ (actually we have to shrink the domain  a bit such 
that the image of $\a$ remains in $U$). It remains to show that $F$ 
extends holomorphically to a \nbd\ of $A$. From 
$$
	F(x)=\bigl(x, \wt f(\a(x))\bigr),
	\quad f(\a(x))= \bigl(\a(x),\wt f(\a(x))\bigr)
$$
we see  that these two points belong to the same leaf of $\cH$. 
Hence $F(x)$ is  the unique point of $N_x$ obtained from $f(\a(x))\in N_{\a(x)}$ 
by a parallel transport along the leaf of $\cH$ through $f(\a(x))$. 
(More precisely, we take the nearest intersection point of the leaf 
with the fiber $N_x$.) Since $\a$ is a biholomorphism close to the 
identity in a \nbd\ of $A$, this gives a well defined \holo\ extension 
of $F$ to a \nbd\ of $A \cup B$ which is transverse to $\cH$. 

Consider now the problem (b). Let $p\in X$ be a critical point of
$\rho$ and assume that $f$ is already defined on $\{\rho\le c\}$ 
for some $c<\rho(p)$ close to $\rho(p)$. The crossing of the critical level 
is localized in a small admissible \nbd\ $U\subset X$ of $p$, except for 
the last step (Subsect.\ 6.4) which uses the noncritical case (a).
We must explain how to extend $f$ smoothly across the handle 
$E\subset U$ attached to $\{\rho\le c\}$ (see Subsect.\ 6.2 
for the details). Using a trivialization $N|_U\simeq U\times \C^q$ 
we have the following situation:

\smallskip
\item{(i)} $f(x)=(x,\wt f(x))$ where $\wt f$ is a \holo\ submersion from 
a \nbd\ of $U\cap \{\rho \le c\}$ to $\C^q$,

\smallskip
\item{(ii)} $\iota \colon N|_U \to  TX|_U = U\times \C^n$ equals
$\iota(x,v)=(x,A_x v)$ $(v\in \C^q)$ where $A_x$ is a complex 
$n\times q$ matrix of rank $q$ depending continuously on $x\in U$, and 

\smallskip
\item{(iii)} $x\to J\wt f(x) \cdotp A_x \in GL_q(\C)$ is 
homotopic to the constant map $x\to I_q$ in a \nbd\ of 
$U\cap \{\rho \le c\}$. 

\smallskip
Note that (iii) is just the condition on $f'\circ\iota$ expressed in local
coordinates. An elementary consequence of (iii) is that the Jacobian 
matrix $J\wt f$ admits a smooth extension across the handle $E\subset \R^k$ 
to  a map $\wt J$ into the space of complex $n\times q$ matrices of rank $q$ 
such that $x\to \wt J(x)\cdotp A_x \in GL_q(\C)$ 
remains homotopic to the constant map on the set 
$(\{\rho\le c\} \cup E) \cap U$. 

Let $D$ be a domain in $\R^n=\R^n+i0 \subset \C^n$ containing 
the handle $E$ as in Lemma 6.4. Let $\Omega$ denote the differential relation 
of order one whose holonomic sections are smooth maps $h\colon  D\to \C^q$ 
whose Jacobian satisfies the condition $Jh(x) \cdotp A_x\in GL_q(\C)$. 
We see as in Lemma 6.5 above that $\Omega$ is ample in the coordinate directions.
Hence Gromov's convex integration lemma from [Gro3, 2.4.1.] (or Sect.\ 18.2 in [EM])
gives a smooth extension of $\wt f$ across the handle $E$ such that 
$x\to J\wt f(x)\cdotp A_x \in GL_q(\C)$ is homotopic to identity on 
$(\{\rho\le c\}\cup E)\cap U$, thereby insuring that the extended section
$f(x)=(x,\wt f(x))$ is transverse to $\cH$ also over $E$ and 
$f'\circ\iota$ remain homotopic to the identity on $N$. (See Lemma 6.5
for the details.) The remaining steps of the proof are the same 
as for submersions. Theorem 7.1 is proved.

\proclaim COROLLARY 7.2: 
Let $V$ be a closed complex submanifold in a Stein manifold $X$.
If the tangent bundle $TX$ admits a complex vector subbundle $N$ with locally constant 
transition functions such that $TX|_V= TV\oplus N|_V$ then $V$ is a union 
of leaves in a nonsingular holomorphic foliation of $X$.

\demo Proof: Let $\cH$ be a foliation of $N$ as in the proof of 
Theorem 7.1. By the Docquier-Grauert theorem [DG] there are an 
open \nbd\ $U\subset X$ of $V$, a \holo\ retraction 
$\pi\colon U\to V$ and an injective holomorphic map 
$\phi\colon U\to  N|_V$ such that 
$\phi(x)\in N_{\pi(x)}$ for each $x\in U$, and 
$\phi(x)=0_x$ \iff\  $x\in V$. The point $\phi(x)$
corresponds to a unique point $f(x)\in \pi^*(N|_V)_x$
via the pull-back map $\pi^*$. Shrinking $U$ if necessary we 
have $N|_U \simeq \pi^*(N|_V)$. Using this identification 
we see that $f\colon U\to N|_U$ is a \holo\ section which intersects
the zero section of $N$ transversely along $V$. 
Shrinking $U$ again we conclude that $f$ is transverse 
to $\cH$ and $V$ is a leaf of the associated foliation $\cF_f$
of $U$.  It remains to find a global transverse section 
$\wt f\colon X\to N$ which agrees with $f$ to second order along $V$. 
This is done as in the proof of Theorem 2.5 (Subsect.\ 6.5),
with the modifications explained above.

\medskip 
{\it Remark.} A closed connected complex submanifold $V$ in a Stein 
manifold $X$ is a leaf in a nonsingular \holo\ foliation defined in 
an open \nbd\ of $V$ \iff\ the normal bundle of $V$ in $X$ admits 
locally constant transition functions. The proof is essentially 
the same as for smooth foliations:
the direct part is due to Ehresmann (see e.g.\ [God, p.\ 5]);
for the converse part we transfer the above foliation $\cH$ of the
normal bundle $N$ to a \nbd\ of $V$ in $X$ by the Docquier-Grauert theorem [DG].
Corollary 7.2 gives a sufficient condition for the 
existence of a {\it global foliation} of $X$ with the same property.

\medskip
{\it Acknowledgements.}
I wish to thank the colleagues who have contributed to this 
work through useful discussions: B.\ Drinovec, P.\ Ebenfelt, P.\ Heinzner, 
L.\ Lempert, E.\ L\o w,  T.\ Ohsawa, N.\ \O vrelid, J.\ Prezelj,
M.\ Slapar, and J.\ Winkelmann. I especially thank N.\ \O vrelid 
who carefully read the paper and kindly pointed out to me a few mistakes.
I sincerely thank R.\ Narasimhan for having initiated my interest in these questions, 
and R.\ Gunning for his interest in a preliminary version of this work 
presented at the Gunning/Kohn conference in Princeton on September 20, 2002. 
I am grateful to the referee for the constructive criticism and for very 
helpful remarks and suggestions. Finally, I thank my colleagues in the 
Mathematics department of the University of Ljub\-ljana 
for providing a pleasant and stimulating environment.
This research has been supported in part by a grant from the 
Ministry of Science and Education of the Republic of Slovenia. 

%
%
%
%
\medskip
\def\same{{--------\ }}

\beginsection References

\item{[A]} 
ANDERS\'EN, E., 
Volume-preserving automorphisms of $\C^n$.
{\it Complex Variables}, 14  (1990), 223-235.

\item{[AL]} 
\same, LEMPERT, L.,
On the group of holomorphic automorphisms of $\C^n$.
{\it Inventiones Math.}, 110  (1992), 371--388.

\item{[AF]} 
ANDREOTTI, A., FRANKEL, T.,
The Lefschetz theorem on hyperplane sections.
{\it Ann.\ Math.}, 69 (1959), 713--717.

\item{[BN]} 
BELL, S., NARASIMHAN, R.,
Proper holomorphic mappings of complex spaces, in 
{\it Several complex variables}, Vol.\ VI, pp.\ 1--38. 
Encyclopaedia Math.\ Sci., 69. Springer, Berlin, 1990.

\item{[DG]} 
DOCQUIER, F., GRAUERT, H.,
Levisches Problem und Rungescher Satz f\"ur 
Teilgebiete Steinscher Mannigfaltigkeiten. 
{\it Math.\ Ann.}, 140 (1960), 94--123.

\item{[EM]}
ELIASHBERG, Y., MISHACHEV, N., 
{\it Introduction to the $h$-principle.} 
Graduate Studies in Math., 48.
Amer.\ Math.\ Soc., Providence, RI, 2002.

\item{[Fe]} 
FEIT, S.,
$k$-mersions of manifolds.
{\it Acta Math.}, 122 (1969), 173--195.

\item{[Fo1]} 
FORSTER, O.,
Some remarks on parallelizable Stein manifolds. 
{\it Bull.\ Amer.\ Math.\ Soc.}, 73 (1967), 712--716.

\item{[Fo2]}
\same, 
Plongements des vari\'et\'es de Stein.  
{\it Comment.\ Math.\ Helv.}, 45 (1970), 170--184.

\item{[F1]} 
FORSTNERI\v C, F., 
Approximation by automorphisms on smooth submanifolds of $\C^n$.
{\it Math.\ Ann.}, 300 (1994), 719--738.

\item{[F2]} 
\same, 
Interpolation by holomorphic automorphisms and embeddings in $\C^n$.
{\it J.\ Geom.\ Anal.}, 9 (1999), 93-118.

\item{[F3]}  \same 
The homotopy principle in complex analysis: A survey.
(Explorations in Complex and Riemannian Geometry: A Volume 
dedicated to Robert E.\ Greene, 
Ed.\ J.\ Bland, K.-T.\ Kim, S.\ G.\ Krantz).
{\it Contemporary Math.}, American Mathematical Society, Providence, 2003.

\item{[FL]} 
\same, L\O W, E.,
Global holomorphic equivalence of smooth submanifolds
in $\C^n$. 
{\it Indiana Univ.\ Math.\ J.}, 46 (1997), 133--153.

\item{[FLO]}
\same, L\O W, E., \O VRELID, N.,
Solving the $d$- and $\overline\partial$-equa\-tions in thin tubes 
and applications to mappings. 
{\it Michigan Math.\ J.}, 49 (2001), 369--416. 

\item{[FP1]} \same, PREZELJ, J., 
Oka's principle for holomorphic fiber bundles with sprays.
{\it Math.\ Ann.}, 317 (2000), 117-154.

\item{[FP2]} 
\same,  
Oka's principle for holomorphic submersions with sprays.
{\it Math.\ Ann.}, 322  (2002), 633-666.

\item{[FP3]} 
\same, 
Extending holomorphic sections from complex subvarieties.
{\it Math.\ Z.}, 236  (2001), 43--68.

\item{[FR]} 
\same, ROSAY, J.-P.,
Approximation of biholomorphic mappings by automorphisms of $\C^n$.
{\it Invent.\ Math.}, 112  (1993), 323--349.
Erratum in {\it Invent.\ Math.}, 118  (1994), 573--574.

\item{[God]} 
GODBILLON, C.,
{\it Feuilletages, \'etudes g\'eom\'etriques}.
Birkh\"auser, Basel-Boston-Berlin, 1991.

\item{[GG]}
GOLUBITSKY, M., GUILLEMIN, V., 
{\it Stable Mappings and their Singularities.} 
Graduate Texts in Mathematics, 14. 
Springer-Verlag, New York-Heidelberg, 1973.
%
%

\item{[Gra]} 
GRAUERT, H., 
Analytische Faserungen \"uber holomorph-vollst\"andigen \break
R\"aumen.
{\it Math.\ Ann.}, 135 (1958), 263--273.

\item{[Gro1]} 
GROMOV, M.,
Stable maps of foliations into manifolds.
{\it Izv.\ Akad.\ Nauk, S.S.S.R.}, 33 (1969), 707--734.

\item{[Gro2]} 
\same,
Convex integration of differential relations, I. (Russian) 
{\it Izv.\ Akad.\ Nauk SSSR Ser. Mat.}, 37 (1973), 329--343. 
English translation in {\it Math.\ USSR--Izv.}, 7 (1973), 329--343.

\item{[Gro3]} 
\same,  
{\it Partial Differential Relations.}
Ergebnisse der Mathematik und ihrer Grenzgebiete (3), 9.
Springer, Berlin--New York, 1986.

\item{[Gro4]} 
\same, 
Oka's principle for holomorphic sections of elliptic bundles.
{\it J.\ Amer.\ Math.\ Soc.}, 2 (1989), 851-897.

\item{[GN]} 
GUNNING, R.\ C., NARASIMHAN, R., 
Immersion of open Riemann surfaces. 
{\it Math.\ Ann.}, 174  (1967), 103--108. 

\item{[GR]} 
\same, ROSSI, H., 
{\it Analytic Functions of Several Complex Variables.}
Pren\-tice--Hall, Englewood Cliffs, 1965.

\item{[Ha1]} 
HAEFLIGER, A., Vari\'et\'es feuillet\'es. 
{\it Ann.\ Scuola Norm.\ Sup.\ Pisa} (3), 16 (1962),  367--397. 

\item{[Ha2]} 
\same  
Lectures on the theorem of Gromov, in  
{\it Proceedings of Liverpool Singularities Symposium}, 
Vol.\ II (1969/1970), pp.\ 128--141. 
Lecture Notes in Math., 209. Springer, Berlin, 1971. 

\item{[HW1]} 
HARVEY, F.\ R., WELLS, R.\ O., Jr.,
Holomorphic approximation and hyperfunction theory on a
$\cC^1$ totally real submanifold of a complex manifold.
{\it Math.\ Ann.}, 197 (1972), 287--318.

\item{[HW2]} 
\same, 
Zero sets of non-negative strictly plurisubharmonic functions.
\break
{\it Math.\ Ann.}, 201 (1973), 165--170.

\item{[HL1]} 
HENKIN, G.\ M., LEITERER, J.,
{\it Theory of Functions on Complex Manifolds.}
Akademie-Verlag, Berlin, 1984.

\item{[HL2]} 
\same,
{\it Andreotti-Grauert Theory by Integral Formulas.}
Progress in Math., 74, Birkh\"auser, Boston, 1988.

\item{[HL3]} 
\same,  
The Oka-Grauert principle without induction over the basis dimension.
{\it Math.\ Ann.}, 311 (1998), 71--93.

\item{[Hi1]} HIRSCH, M.,  
Immersions of manifolds.
{\it Trans.\ Amer.\ Math.\ Soc.}, 93 (1959), 242--276.

\item{[Hi2]} \same,
On embedding differential manifolds into Euclidean space.
{\it Ann.\ Math.}, 73 (1961), 566--571.

\item{[H\"o1]} 
H\"ORMANDER, L., 
$L\sp{2}$ estimates and existence theorems for the $\bar \partial$ operator. 
{\it Acta Math.}, 113 (1965), 89--152. 

\item{[H\"o2]} \same, 
{\it An Introduction to Complex Analysis in Several Variables}. 
Third ed. 
North Holland, Amsterdam, 1990.

\item{[HWe]} 
\same, WERMER, J.,
Uniform approximations on compact sets in $\C^n$.
{\it Math.\ Scand.}, 23 (1968), 5--21.

\item{[Ko]} 
KOLMOGOROV, A.,
On representation of continuous functions of several variables
by superpositions of continuous functions of fewer variables.
{\it Dokl.\ Akad.\ Nauk SSSR}, 108 (1956), 179--182.

\item{[Ku]} 
KUTZSCHEBAUCH, F., 
Anders\'en-Lempert theory with parameters.
Pre\-print, 2002.

\item{[LV]} 
LEHTO, O., VIRTANEN, K.\ I.,
{\it Quasiconformal Mappings in the Plane}. Second ed.
Grundlehren der math.\ Wiss., 126. 
Springer, New York--Heidelberg, 1973.

\item{[MS]}
MILNOR, J.\ W., STASHEFF, J.\ D.,
{\it Charactersitic Classes}.
Annals of Math.\ Studies, 76.
Princeton Univ.\ Press, Princeton, 1974.

\item{[N]} 
NISHIMURA, Y., 
Examples of analytic immersions of two-dimensional Stein 
manifolds into $\C^2$. 
{\it Math.\ Japon.}, 26 (1981), 81--83.

\item{[Pf]} 
PFLUGER, A., 
\"Uber die Konstruktion Riemannscher Fl\"achen durch Verheftung.
{\it J.\ Indian Math.\ Soc.} (N.S.), 24 (1961), 401--412.

\item{[Ph1]} 
PHILLIPS, A.,
Submersions of open manifolds.
{\it Topology}, 6 (1967), 170--206.

\item{[Ph2]} 
\same,
Foliations on open manifolds, I. 
{\it Comm.\ Math.\ Helv.}, 43 (1968), 204--211.

\item{[Ph3]} 
\same,
Foliations on open manifolds, II.
{\it Comm.\ Math.\ Helv.} 44 (1969), 367--370.

\item{[Ph4]} 
\same,
Smooth maps transverse to a foliation. 
{\it Bull.\ Amer.\ Math.\ Soc.}, 76 (1970), 792--797.

\item{[Ph5]}  
\same, 
Smooth maps of constant rank. 
{\it Bull.\ Amer.\ Math.\ Soc.}, 80 (1974), 513--517.

\item{[Ra]} RAMSPOTT, K.\ J.,
Stetige und holomorphe Schnitte in B\"undeln mit homogener Faser.
{\it Math.\ Z.}, 89 (1965), 234--246.
 
\item{[RS]}  
RANGE, R.\ M., SIU, Y.\ T., 
$\cC^k$ approximation by holomorphic functions and $\bar \partial$-closed 
forms on $\cC^k$ submanifolds of a complex manifold.
{\it Math.\ Ann.}, 210 (1974), 105--122.

\item{[Ro]} 
ROSAY, J.-P.,
A counterexample related to Hartog's phenomenon
(a question by E.\ Chirka). 
{\it Michigan Math.\ J.}, 45 (1998), 529--535.

\item{[Sm]}
SMALE, S.,
The classification of immersions of spheres in Euclidean spa\-ces. 
{\it Ann.\ of Math.} (2), 69 (1959), 327--344.

\item{[Sp]} 
SPRING, D., 
{\it Convex Integration Theory.}   
(Solutions to the $h$-principle in geometry and topology).
Monographs in Mathematics, 92.
Birk\-h\"auser, Basel, 1998.

\item{[Ste]} 
STEIN, K.,  
Analytische Funktionen mehrerer komplexer Ver\"anderlichen zu 
vorge\-ge\-benen Periodizit\"atsmoduln und das zweite Cousinsche Pro\-blem.
{\it Math.\ Ann.}, 123 (1951), 201--222.

\item{[Sto]}
STOUT, E.\ L.,
{\it The Theory of Uniform Algebras.} 
Bogden\&Quigley, Inc., Tarrytown-on-Hudson, New York, 1971. 

\item{[Tho]}
TH\'OM, R., 
Un lemme sur les applications diff\'erentiables.  
{\it Bol.\ Soc.\ Mat.\ Mexicana}, (2) 1 (1956), 59--71.

\item{[Th1]} 
THURSTON, W.,  
The theory of foliations of codimension greater than one. 
{\it Comm.\ Math.\ Helv.}, 49 (1974), 214--231.

\item{[Th2]} 
\same, 
Existence of codimension one foliations.
{\it Ann.\ Math.}, 104 (1976), 249--268.

\item{[V]} 
VAROLIN, D., 
The density property for complex manifolds and geometric structures.
{\it J.\ Geom.\ Anal.}, 11 (2001), 135--160.

\bye
 
\bigskip\medskip
{\it Address:}
Institute of Mathematics, Physics and Mechanics,
University of Ljubljana, Jadranska 19, 1000 Ljubljana, Slovenia

\medskip
{\it E-mail:} franc.forstneric@fmf.uni-lj.si

\bye